\documentclass[reqno,11pt]{article}
\usepackage{a4wide,xcolor,eucal,enumerate,mathrsfs}
\usepackage[normalem]{ulem}
\usepackage[pdfborder={0 0 0}]{hyperref}  

\usepackage{amsmath,amssymb,epsfig,amsthm}
\usepackage[latin1]{inputenc}

\numberwithin{equation}{section}

\newtheorem {theorem}{Theorem}[section]
\newtheorem {Corollary}[theorem]{Corollary}
\newtheorem {Definition}[theorem]{Definition}
\newtheorem {Lemma}[theorem]{Lemma}

\newtheorem {Remark}[theorem]{Remark}
\newcommand{\R}{{\mathbb R}}
\newcommand{\N}{{\mathbb N}}
\newcommand{\g}{{{\rm g}}}

\title{Embedded surfaces of arbitrary genus minimizing  the Willmore energy under isoperimetric constraint}

\begin{document}

\author{Laura Gioia Andrea Keller \footnote{Institut f\"ur Numerische und Angewandte Mathematik, Universit\"at M\"unster, 48149 M\"unster, Germany}, 
Andrea Mondino \footnote{Department of Mathematics, ETH Z\"urich, 8092 Z\"urich, Switzerland}, Tristan Rivi\`ere \footnote{Department of Mathematics, ETH Z\"urich, 8092 Z\"urich, Switzerland}
}
\maketitle

\begin{abstract}
The isoperimetric ratio of an embedded surface in $\R^3$ is defined as the ratio of the area of the surface to power three to the squared enclosed volume. \\
The aim of the present work is to study the minimization of the Willmore energy under fixed isoperimetric ratio when
the underlying abstract surface has fixed genus $\g\geq 0$.  \\
The corresponding problem in the case of spherical surfaces, i.e. $\g=0$,  was recently solved by Schygulla (see \cite{Schy}) with different methods.
\end{abstract}

\textit{Keywords:} Willmore functional, Helfrich functional, isoperimetric constraint, tori, higher genus surfaces.


\tableofcontents

\section{Introduction and motivation} 

The Willmore functional of an immersion $\vec{\Phi}$ of an abstract oriented surface $\Sigma^2$  is given by
\begin{displaymath}
W(\vec{\Phi})=\int_{\Sigma^2} \vert \vec{H}\vert^2 \; dvol_g
\end{displaymath}
 where $\vec{H}$ is the mean curvature of the immersion $\vec{\Phi}$ and $dvol_g$ the induced volume form.\\
A surface which is a critical point of this functional is called Willmore surface. These surfaces have been introduced by Wilhelm Blaschke in the beginning of the XX's century in the framework of ``conformal geometry'' and were called at the time ``conformal minimal surfaces'' (see \cite{Bla}). He observed that the Willmore energy has the property to be invariant under conformal transformations of the ambient spaces, the transformations which infinitesimally preserve angles. He introduced then the theory of Willmore surfaces as being the ``natural'' merging between minimal surface theory and conformal invariance. 
 \\

 Probably because of the richness of the symmetries preserving the class of Willmore surfaces and because of the simplicity and the universality of its definition, the Willmore functional shows up in various fields of sciences and technology. It appears for instance in biology in the study of lipid bilayer cell membranes
under the name ``Helfrich energy''  (as we will see in more detail later in the introduction),   in general relativity as being the main term in the so called ``Hawking Mass'', in string theory in high energy physics it appears in the definition of the
Polyakov extrinsic action, in elasticity theory as  free energy of the non-linear plate Birkhoff theory, in optics and lens design, \dots etc.
 
 The theory of Willmore surfaces, named after the person who revisited the theory in an important work of the mid-sixties ( see \cite{Wil} a\cite{Wil2}), has flourished
 in the last decades.\\

Over the last twenty years, existence of minimizers of the Willmore energy under various constraints has been
obtained: in the class of smooth immersions for a given abstract surface $\Sigma^2$
(see \cite{Si} and  \cite{BK}),  within a fixed conformal class (see \cite{KL} and \cite{KS6}) or more recently
 in the class of smooth embeddings of the sphere under fixed isoperimetric ratio (see \cite{Schy}) 
 \begin{equation}\label{eq:defIso}
iso(\vec{\Phi})=\frac{(Area(\vec{\Phi}))^3}{\Big(Vol\Big( \; \overline{\vec{\Phi}} \;\Big)\Big)^2}
\end{equation}
where $Vol\Big( \; \overline{\vec{\Phi}} \; \Big)$ denotes the volume enclosed by $\vec{\Phi}(\Sigma^2)$, i.e. the volume of the bounded
connected component of $\mathbb{R}^3 \backslash \vec{\Phi}(\Sigma^2)$.
\\

The aim of the present work is to study the minimization of the Willmore energy under fixed isoperimetric ratio when
the underlying abstract surface has fixed genus $\g\geq 0$.\\

Beside the obvious geometric intrinsic interest such a minimization under  isoperimetric and genus constraint could have, a motivation to study this problem comes from the modelization of the free energy of elastic lipid bilayer membranes in cell biology. 
\\Indeed the Willmore functional is closely related to the Helfrich functional which describes the free energy of a closed lipid bilayer 
\begin{displaymath}
F_{Helfrich}= \int_{\textrm{lipid bilayer}}\Big ( \frac{k_c}{2}(2H+c_0)^2 +\bar{k}K + \lambda \Big ) + p \cdot V
\end{displaymath}
where $k_c$ and $\bar{k}$ denote bending rigidities, $c_0$ stands for the spontaneous curvature, 
$\lambda$ is the surface tension, $K$ and $H$ denote as usual the Gauss curvature 
and the mean curvature, respectively, $p$ denotes the osmotic pressure and $V$ denotes the enclosed volume. \\ 
The shapes of such membranes at equilibrium are then given by the corresponding Euler-Lagrange equation.
If $c_0=0$, $\lambda=0$ and $p=0$ the Willmore functional captures the leading terms in Helfrich's functional (up to a topological constant).
Whereas if these physical constants do not vanish, $\lambda$ and $p$ can be seen as Lagrange
multipliers for area and volume constraints. Thus, thanks to the invariance under rescaling of both the Willmore functional and the isoperimetric ratio, we exactly face the problem of minimizing the Willmore functional under an isoperimetric constraint. \\
In the context of vesicles, imposing a fixed area and a fixed volume has perfect biological meaning:
on one hand, it is observed that at experimental time scales the lipid bilayers exchange only few molecules 
with the ambient and the possible contribution to the elastic energy due to displacements within the
membrane is negligible. Thus, the area of the vesicle can be treated as a fixed one. On the other hand,
a change in volume would be the result of a transfer of liquid into or out of the vesicle. But this would significantly change the osmotic pressure
and thus would lead to an energy change of much bigger scale than the scale of bending energy. \\
At first glimpse one may think that biologically relevant vesicles should always be of spherical shape. 
But in fact also higher genus membranes are observed: for toroidal shapes see \cite{MB} and \cite{Seifert}, for genus two surfaces  see \cite{MB1}, and for higher genuses see \cite{MB2}. Further details  can be found also in \cite{LS}. \\
In addition, often in biology the ratio of area 
(the place where a molecule is produced) to volume (how much of the produced molecule can be stored) is crucial. \\
Moreover, the Helfrich functional as well as the Willmore functional are widely used for modeling biological phenomena, e.g. red blood cells (see e.g. \cite{Mil}), 
folds of the endoplasmic reticulum (see e.g. \cite{Rap}) and morphologies (Cristae junction) of mitochondria (see e.g. \cite{Mitoch}). \\
Thus, also from an applied, biological point of view, 
it is perfectly reasonable to look at the problem we propose to study in this article. 
\\

In order to describe the results of this paper, let us introduce the framework where we are going to work. 
\\

In the aforementioned article of Schygulla (see \cite{Schy}), the main analytical strategy was to study immersions from the point of view of the image $\vec{\Phi}(\Sigma^2)$. 
This strategy, which has been extensively  used also in a series of works by Kuwert and Sch\"atzle (see \cite{KS1}- \cite{KS2}-\cite{KS3}-\cite{KS4}-\cite{KS5}-\cite{KS6}-\cite{KS7}), was introduced by Simon (see \cite{Si}) and it  is also known under the denomination of \emph{Simon's ambient approach}. \\

In a series of papers (see \cite{Riv2}, \cite{Riv5}, \cite{Riv4}), the third author established a new framework for the study of variational problems related to the Willmore functional in which he was favoring the study of the immersion $\vec{\Phi}$ itself instead of
its image. This approach could be called \emph{the parametric approach}. It provides a general framework in which not only the above mentioned existence results (free minimization and minimization in a fixed conformal class) could be extended to the class of Lipschitz immersions with $L^2$-bounded second fundamental form, also called \emph{weak immersions}, but it is also suitable for applying fundamental 
principles of the calculus of variations (such as the mountain pass lemma for instance,  in order to produce saddle type critical points). \\

In the present paper we will adopt this parametric approach. For the reader's convenience, all  the important concepts as e.g. this space of Lipschitz immersions that we will denote ${\mathcal E}_{\Sigma^2}$, are introduced and explained in Section 2. \\

From now on, the underlying abstract surface is  closed   of genus $\g\geq 0$  and  our goal is to look for minimizers of the Willmore energy in the class of Lipschitz immersions with second fundamental form in $L^2$ under the additional constraint that the isoperimetric ratio, defined in \eqref{eq:defIso}, is fixed (by definition, if $\vec{\Phi}\in {\mathcal E}_{\Sigma^2}$ is not an embedding, i.e. it has self intersections, we set $iso(\vec{\Phi}):=+\infty$).
\\
Our first result is an alternative proof, using the parametric approach, of Schygulla's theorem mentioned above (see \cite{Schy}), i.e. existence of smooth embedded spheres minimizing $W$ under isoperimetric constraint. Actually we manage to prove a stronger result, namely our minimization is performed in the larger class of weak immersions, $\mathcal{E}_{{\mathbb S}^2}$, rather than among smooth embeddings.
\\ Before stating it, recall that by the isoperimetric inequality in $\R^3$, for every embedded surface $\vec{\Phi}$ one has $iso(\vec{\Phi})\geq iso(\mathbb{S}^2)=36\pi$ and equality occurs if and only if $\vec{\Phi}$ is actually a round sphere. 

\begin{theorem}[The genus 0 case]\label{thm:Schy}
For every $R\in [36\pi, +\infty)$ there exists a smooth embedded spherical surface, called later on Schygulla sphere and denoted by ${\mathbb S}_{S,R}$,  which minimizes the Willmore functional $W$ among weak immersions ${\mathcal E}_{{\mathbb S}^2}$ having constrained isoperimetric ratio equal to $R$.
\end{theorem}

Since, in virtue of  the theorem above, the genus $0$ case is well understood, from now on we will assume that $\Sigma^2$ is a surface of genus $\g\geq 1$. Before stating our main theorems let us introduce some notation. \\
 Let $\beta^3_{\g}$ denote the infimum (actually it is a minimum thanks to \cite{Si},\cite{BK} and \cite{Riv5}) of the Willmore energy among surfaces of genus $\g\geq 1$ immersed in $\R^3$
\begin{equation}\label{eq:defbeta}
\beta^3_{\g} := \inf \left \{ W(\vec{\Phi}) \; \vert \; \vec{\Phi} \, \textrm{is an immersion of the genus $\g$ closed surface} \right \}\quad .
\end{equation}
Let $\vec{\Phi}_{\g}$ be a smooth minimizer attaining the infimum $\beta^3_{\g}$ (notice that all the compositions of $\vec{\Phi}_{\g}$ with conformal maps of $\R^3$ are still minimizers thanks to the conformal invariance of $W$) and let $\omega^3_{\g}$ be defined by
\begin{equation}\label{eq:defomega}
\omega^3_{\g} := \min \left \{ 4 \pi +  \sum_{i=1}^p (\beta^3_{\g_i} - 4  \pi) \; \vert \; \g=\g_1 + \dots + \g_p \; , \; 1 \leq \g_i < \g \right \} \quad .
\end{equation}
 Now we can state the main results of the paper. The first one guarantees the existence of minimizers of our problem under certain hypothesis and the second one 
discusses when such assumptions are satisfied. 

\begin{theorem}[Isoperimetric-constrained minimizers of $W$ of every genus]\label{MainTheorem}
Let $\Sigma^2_{\g}$ be the abstract closed (i.e. compact without boundary) orientable surface  of genus $\g\geq 1$   and consider the set
\begin{equation}\label{eq:defI}
I_{\g}:=  \left \{ R \in \mathbb{R} \, \Big \vert \, 
\inf_{\begin{subarray}{c}\vec{\Phi} \in \mathcal{E}_{\Sigma_{\g}^2} \\ iso(\vec{\Phi})=R 
\end{subarray}}W(\vec{\Phi})<
\min \left \{ 8 \pi, \omega^3_{\g}, \beta^3_{\g} + W(\mathbb{S}_{S,R}) - 4 \pi \right \}
\right  \}.
\end{equation}
Then for given 
$R \in I_{\g}$ there exists an embedding $\vec{\Phi}$ of $\Sigma^2_{\g}$
into $\mathbb{R}^3$,  with $iso(\vec{\Phi})=R$,
which minimizes the Willmore energy among all Lipschitz immersions of $\Sigma^2_{\g}$ with second fundamental form bounded in $L^2$ 
(i.e. among $\mathcal{E}_{\Sigma^2_{\g}}$) and fixed isoperimetric ratio equal to $R$.
\end{theorem} 

We even know that the minimizer from Theorem \ref{MainTheorem} is a smooth embedding.

\begin{Corollary}[Smoothness of the minimizer] \label{Corollary}
The minimizing $\vec{\Phi}$ from Theorem \ref{MainTheorem} is a smooth embedded genus $g\geq 1$ surface which minimizes the
Willmore energy among smooth genus $\g$ embeddings with given isoperimetric ratio $R$.
\end{Corollary}

As announced above, we will characterize further the set $I_{\g}$ for which there exist minimizers of the Willmore energy under 
the additional constraint of given isoperimetric ratio.

\begin{theorem}[$I_{\g}\neq \emptyset$ is open]\label{Thm2}
$I_{\g} \subset \R$ is a non-empty open set containing $iso(\vec{\Phi}_{\g})$, the isoperimetric ratio of any \emph{free} minimizer of $W$ among smooth genus $\g$ immersed closed surfaces.  
\end{theorem}

\begin{Remark}[$ I_{\g} \supset ( 36 \pi , R_{\g}+\delta) $]\label{rem:Thm2}
\rm{
Notice that from Theorem \ref{Thm2}, $I_{\g}$ contains the whole interval  $(36\pi, iso(\vec{\Phi}_{\g})]$, where $36 \pi= iso({\mathbb S}^2)$. Indeed, since  for \emph{every} minimizer $\vec{\Phi}_{\g}$ of the free minimization problem one has $iso(\vec{\Phi}_{\g}) \in I_{\g}$ and since   the embeddings $\vec{\Xi}\circ \vec{\Phi}_{\g}$ are still minimizers of $W$, where $\vec{\Xi}$ is any conformal transformation (i.e. a M\"obius map) of $\R^3$, it follows that $iso(\vec{\Xi}\circ \vec{\Phi}_{\g})\in I_{\g}$ for every $\vec{\Xi}$ . 
\\
Let us now consider a special smooth 1-parameter family $\{\vec{\Xi}_r\}_{r\in (0,+\infty)}$ of M\"obius maps given by the inversion with respect to a sphere of unit radius and center $p(r)$, where $p:(0,\infty)\to \R^3$ is a smooth curve  such that, for every $r>0$, $p(r)$ is  at distance $\frac{1}{r}$ from the fixed minimizer ${\Phi}_{\g}(\Sigma^2_{\g}).$ 
\\
From this construction, it is not difficult to check that $iso(\vec{\Xi}_r\circ \vec{\Phi}_{\g})$ varies smoothly on $r$ and that
$$\lim_{r\to 0^+}iso(\vec{\Xi}_r\circ \vec{\Phi}_{\g})=iso(\vec{\Phi}_{\g})\quad \text{and} \quad    \lim_{r\to +\infty} iso(\vec{\Xi}_r\circ \vec{\Phi}_{\g})=iso({\mathbb S}^2)=36\pi.$$
A continuity argument then shows that $I_{\g}\supset (36\pi, iso(\vec{\Phi}_{\g})]$. 
\\On the other hand, observe that for $R \in (36\pi, iso(\vec{\Phi}_{\g})]$, the isoperimetric-constrained minimizer produced by Theorem \ref{MainTheorem} is exactly the free minimizer $\Xi_{r(R)}\circ\vec{\Phi}_{\g}$, for a suitable $r=r(R)$.
Nevertheless, called 
$$R_{\g}:=\max_{\vec{\Xi}\in Moeb(\R^3)} iso(\vec{\Xi}\circ \vec{\Phi}_{\g})$$
(notice that the  maximum is attained since, as explained above, when the  M\"obius map diverges the isoperimetric ratio converges to the minimum value, given by  $36\pi=iso(\mathbb{S}^2$)),  the point of Theorem \ref{Thm2} is that $I_{\g}$ is \emph{open}, so it contains an interval of the type $(36\pi, R_{\g}+\delta)$ for some $\delta>0$; and in the interval $(R_{\g}, R_{\g}+\delta) $ the constrained minimizer is a new surface (i.e. not free Willmore) which is interesting to investigate. Let us also stress that we expect that the conditions used to define $I_{\g}$ are satisfied by a larger class of isoperimetric ratios, for more details see Remark \ref{rem:Assumptions}.  
} 
 \hfill$\Box$
\end{Remark}

Since the genus-one-case, i.e. $\Sigma^2=\mathbb{T}^2$ is the 2-d torus, is particularly important for  applications, let us discuss it in more detail. Thanks to the recent proof of the Willmore conjecture by Marques and Neves (see \cite{MN}), we know that $\beta^3_1=2\pi^2$, the set of  minimizers is made by  the Clifford torus ${\mathbb T}^2_{Clifford}$ (i.e. the torus of revolution with radii ratio $1:\sqrt{2}$) and its images under conformal mappings of $\R^3$. Moreover, by definition, $\omega^3_1=+\infty$. We summarize Theorem \ref{MainTheorem}, Corollary \ref{Corollary} and Theorem \ref{Thm2} for the genus-one-case in the following theorem.

\begin{theorem}[Isoperimetric-constrained minimizers of $W$ among tori]\label{thm:MinTorus}
Let ${\mathbb{T}}^2$ be the abstract $2$-dimensional torus and consider the set
\begin{equation}\nonumber
I_1:=  \left \{ R \in \mathbb{R} \, \Big \vert \, 
\inf_{\begin{subarray}{c}\vec{\Phi} \in \mathcal{E}_{{\mathbb T}^2} \\ iso(\vec{\Phi})=R 
\end{subarray}}W(\vec{\Phi})<
\min \left \{ 8 \pi, 2 \pi^2 + W(\mathbb{S}_{S,R}) - 4 \pi \right \}
\right  \}.
\end{equation}
Then 
\begin{itemize}
\item[i)] $I_1 \subset \R$ is a non empty open set satisfying 
$$I_1\supset\left(36\pi, 16 \sqrt{2}\pi^2\right] \quad ,$$
where the numbers above come from the fact that $iso({\mathbb S}^2)=36\pi$ and $iso({\mathbb T}^2_{Clifford})=16\sqrt{2} \pi^2$.
\item[ii)] For given $R \in I_1$ there exists an embedding $\vec{\Phi}$  of  $\,{\mathbb T}^2$ into $\mathbb{R}^3$,  with $iso(\vec{\Phi})=R$,
which minimizes the Willmore energy among all Lipschitz immersions of ${\mathbb T}^2$  with second fundamental form bounded in $L^2$ 
(i.e. among $\mathcal{E}_{{\mathbb T}^2}$) and fixed isoperimetric ratio equal to $R$.
\item[iii)] The minimizing $\vec{\Phi}$ given in $ii)$ is a smooth embedded torus which minimizes the
Willmore energy among smooth embedded tori with isoperimetric ratio $R$.
\end{itemize}
\end{theorem}

\begin{Remark}[About the assumptions of Theorem \ref{MainTheorem}]\label{rem:Assumptions}
\rm{\begin{itemize}

\item[i)] \emph{The $8\pi$-bound:} The natural framework for studying the isoperimetric constraint is given by embedded surfaces (i.e. surfaces without self-intersection) since for this class it is clear what the enclosed volume is. A celebrated inequality of Li and Yau (see Theorem \ref{LiYau}) states that if $\vec{\Phi}$ has self intersections  then $W(\vec{\Phi})\geq 8 \pi$. Therefore, since we want to handle embedded surfaces, it is natural to work under the assumption that the Willmore energy is bounded from above by $8 \pi$. 

\item [ii)] \emph{The $\omega^3_{\g}$-bound:} This assumption is technical and it is used, together with the $8 \pi$-bound, to ensure that the conformal structures of the minimizing sequence do not degenerate in the Moduli space, see Theorem  \ref{ConfClass} in the Appendix. Notice that in the genus-one-case this condition is not needed. An interesting open problem is to rule out this condition from Theorem  \ref{MainTheorem} and study also the possible degenerations.

\item[iii)]\emph{The $\beta^3_{\g}+ W(\mathbb{S}_{S,R}) - 4 \pi$-bound:} Notice that the condition 
\begin{equation}\label{eq:srictIne}
\inf_{\begin{subarray}{c}\vec{\Phi} \in \mathcal{E}_{\Sigma_{\g}^2} \\ iso(\vec{\Phi})=R 
\end{subarray}}W(\vec{\Phi})<  \beta^3_{\g} + W(\mathbb{S}_{S,R}) - 4 \pi 
\end{equation}
is that the \emph{strict} inequality holds.  In fact, for every given isoperimetric ratio different from the one of the round sphere the weak inequality
$$\inf_{\begin{subarray}{c}\vec{\Phi} \in \mathcal{E}_{\Sigma_{\g}^2} \\ iso(\vec{\Phi})=R 
\end{subarray}}W(\vec{\Phi}) \leq  \beta^3_{\g} + W(\mathbb{S}_{S,R}) - 4 \pi $$
holds. This can be seen as follows.
\\ Let $r$ and $\varepsilon > 0$ be given. Then by gluing a Schygulla sphere corresponding 
to an isoperimetric ratio $\tilde{r}=r - \delta$, for some small $\delta>0$, and an inverted minimizer $\vec{\Phi}_{\g}$, i.e. the image of  $\vec{\Phi}_{\g}$ ( recall that  $\vec{\Phi}_{\g}$ is a minimizer of the Willmore energy among genus $\g$ immersed surfaces in $\R^3$, i.e. $W(\vec{\Phi}_{\g})=\beta^3_{\g}$)  under a sphere inversion based at a point of $\vec{\Phi}_{\g}(\Sigma^2_{\g})$ ,
one obtains a surface whose isoperimetric ratio equals the given $r$. 
This can be achieved by appropriately choosing $\delta$ and the way how the two 
parts are glued together. Moreover, the Willmore energy of this new surface can be bounded 
from above by $\beta^3_{\g} + W(\mathbb{S}_{S,r}) - 4 \pi + \varepsilon$. It is an interesting open problem whether or not the strict inequality \eqref{eq:srictIne} is actually always satisfied.
\end{itemize}
} \hfill$\Box$
\end{Remark}

Now we briefly outline the  strategy to prove Theorem \ref{MainTheorem}. As mentioned above, following the approach of  \cite{Riv5}, we will directly work with the immersions $\vec{\Phi}_k$ and not with the immersed surfaces, 
i.e. the images $\vec{\Phi}_k(\Sigma^2) \subset \mathbb{R}^3$. \\
At first glimpse, it may appear a drawback to work with the immersions itself due to the invariance of the problem under the action of the  noncompact  group of  diffeomorphisms of $\Sigma^2$ but, as we will see, 
this difficulty can be handled by choosing an appropriate gauge, namely the \emph{Coulomb gauge}. 
Once such a Coulomb gauge is at hand, one can construct \emph{conformal coordinates}. Together with suitable \emph{estimates for the conformal factors} 
this leads to a setting in which we have very powerful analytical tools at hand in order to solve our problem. 
For example, we can use the fact that the equation for the Willmore surfaces can be reformulated as a conservation law, which enables us to convert the initial 
supercritical problem into a \emph{critical} one.\\

More precisely, fixed $g\geq 1$, we consider a minimizing sequence $\left \{ \vec{\Phi}_k \right \}_{k \in \mathbb{N}}$ in $\mathcal{E}_{\Sigma_{\g}^2}$ for the Willmore energy under the constraint that
the isoperimetric ratio $iso(\vec{\Phi}_k)$ is fixed and equals $R$, for some $R\in I_{\g}$ where $I_{\g}$ is defined in \eqref{eq:defI}.
Notice that, since our abstract surface is a surface of genus $\g\geq 1$, the reference metric of constant scalar curvature is flat or hyperbolic.  Moreover up to a change of coordinates we can assume  that the $\vec{\Phi}_k$s are conformal immersions (see Theorem \ref{conf}) and, since by assumption $W(\vec{\Phi}_k)\leq \min\{8 \pi,\omega^3_{\g}\}-\delta$ ($k$ big enough, for some $\delta>0$),  the conformal classes do not degenerate (see Theorem \ref{ConfClass} in the Appendix). \\

The proof then splits into two cases: either the conformal factors of the immersions $\vec{\Phi}_k$ remain bounded or they diverge to $- \infty$. \\
In the first case, we proceed as follows: in a first step, we show that $\vec{\Phi}_k$ converges weakly in $W^{2,2}$ possibly away from finitely many points of energy concentration. Then, in a second step, we show a point removability result for these possible points of energy concentration. Once we have this, together with the exclusion of possible bubbling, we can conclude convergence in $\mathcal{E}_{\Sigma_{\g}^2} $ and finally show smoothness of the minimizer.
\\
In the second case,  our analysis  reveals that the assumption of diverging conformal
factors yields a dichotomy between one part of the surface carrying the topological information
and another part carrying the isoperimetric ratio; at this point we will perform a ``cut-and-replace'' argument (which recently appeared in a paper of the second and third authors \cite{MR}) which will lead to a contradiction to the hypothesis \eqref{eq:defI} on $I_{\g}$. Thus, the case of diverging conformal factors is excluded and the proof is complete.

\textbf{Acknowledgment}
Part of this work was done when the first two authors were visiting ETH Z\"urich. They would like to thank the FIM (Forschungsinstitut f\"ur Mathematik) and  ETH Z\"urich for the hospitality and the excellent working atmosphere.

\section{Setup and preliminaries}

In this section we will recall the necessary definitions and we will collect and present the most important result which will play a crucial role in our analysis.

\subsection{Definition of the Willmore functional, its fundamental properties and the Willmore surface equation}

First of all, let us recall the definition of the Willmore functional, also called Willmore energy:\\
Let $\Sigma^2$ be an abstract oriented closed surface of genus $\g \geq 0$ (sometimes we will  write $\Sigma^2_{\g}$ to stress the genus), and let $\vec{\Phi}$ be a $C^2$-immersion of the surface $\Sigma^2$ into $\mathbb{R}^3$. 
We denote with $g$ metric induced by the immersion $\vec{\Phi}$, i.e. $g_{ij}:=g_{\R^3}\left( \partial_{x^i} \vec{\Phi} , \partial_{x^j} \vec{\Phi} \right)$; 
in a more geometric language we say that  $g=\vec{\Phi}^* g_{\mathbb{R}^3}$ is  the pullback of the canonical metric on $\mathbb{R}^3$  on $T\Sigma^2$ via the immersion $\vec{\Phi}$, $g$ is also called first fundamental form. There is also a second fundamental form associated to the immersion $\vec{\Phi}$ and it is given by the following map
\begin{eqnarray*}
\vec{\mathbb{I}}_p : T_p\Sigma^2 \times T_p\Sigma^2 &\rightarrow& (\vec{\Phi}_*T_p \Sigma^2 )^{\perp} \\
(X,Y) &\mapsto& \vec{\mathbb{I}}_p(X,Y) := \pi_{\vec{n}}(d^2 \vec{\Phi}(X,Y))=\bar{\nabla}_{\vec{Y}}\vec{X}- \nabla _Y X
\end{eqnarray*}
where $\vec{Z}=d\vec{\Phi}\cdot Z$, $V^{\perp}$ is the orthogonal complement of the subspace $V$ of $\R^3$, $\bar{\nabla}$ is the Levi-Civita connection in $\R^3$ for $g_{\mathbb{R}^3}$ (i.e. the usual directional derivative in $\R^3$)   and $\nabla$ is the Levi-Civita connection on $T\Sigma^2$ induced by $g$. \\

Recall that $\vec{\Phi}$ is said \emph{conformal immersion} if  $\partial_x \vert \vec{\Phi} \vert \equiv \vert \partial_y\vec{\Phi} \vert$ and $g_{\R^3}(\partial_x\vec{\Phi}, \partial_y \vec{\Phi})\equiv 0$. In this case $e^{\lambda}=\vert \vec{\Phi}_x \vert = \vert \vec{\Phi}_y \vert$ is called the \emph{conformal factor} of  $\vec{\Phi}$. \\

Now, we can define the Willmore functional

\begin{Definition} [Willmore Functional]
Let $\Sigma^2$,  $\vec{\Phi}$  and $g$ be as above. Then the \emph{Willmore functional} of the immersion  $\vec{\Phi}$ is given by the following expression
\begin{displaymath}
W(\vec{\Phi})=\int_{\Sigma^2} \vert \vec{H}\vert^2 \; dvol_g \quad ,
\end{displaymath}
where $\vec{H}:=\frac{1}{2}tr(g^{-1}\vec{\mathbb{I}})=\frac{1}{2}\sum_{i,j=1}^2 g^{ij}\vec{\mathbb{I}} 
(\partial_{x_i},\partial_{x_j})$ is the \emph{mean curvature}, and $dvol_g$ is the volume form associated to $g$.
\end{Definition}
Note that we can write $\vec{H}= H \vec{n}$, where $\vec{n}$ is the  unit vector orthogonal to $\vec{\Phi}_*(T\Sigma^2)$ in $\R^3$ given by 
$$\vec{n}:=\frac{\partial_{x^i} \vec{\Phi} \times \partial_{x^j} \vec{\Phi} }{|\partial_{x^i} \vec{\Phi} \times \partial_{x^j} \vec{\Phi}|}, $$
where $\vec{v_1}\times \vec{v_2}$ is the usual  cross product of vectors in $\R^3$.
\\

Recall that the above Willmore functional is \emph{conformally invariant} (see \cite{C} and \cite{RivNotes}), i.e. it is invariant under isometries, dilations and sphere inversions of the ambient space $\R^3$. \\

Moreover,  thanks to the Gauss-Bonnet theorem, we have 
\begin{eqnarray} \label{WillmoreFunctional}
W(\vec{\Phi}) &=& \frac{1}{4}\int_{\Sigma^2} \vert \Delta_g \vec{\Phi} \vert^2  dvol_g = \frac{1}{4}\int_{\Sigma^2} \vert d\vec{n}\vert^2_g dvol_g + \pi \chi(\Sigma^2)  \\
              &=&  \frac{1}{4}\int_{\Sigma^2} \vert \vec{\mathbb{I}} \vert^2 dvol_g + \pi \chi(\Sigma^2) =  \frac{1}{4} \int_{\Sigma^2} (\kappa_1 - \kappa_2)^2 dvol_g + 2\pi \chi(\Sigma^2) \nonumber \quad, 
\end{eqnarray}
where $\chi(\Sigma^2)=2g-2$ denotes the Euler Characteristic of the surface $\Sigma^2$, $\kappa_i$ are the principal curvatures, and 
$\Delta_g$ denotes the intrinsic negative Laplace-Beltrami operator. \\

Next, recall the notion of Willmore immersion (or Willmore surface).
\begin{Definition}[Willmore immersion]
Let $\Sigma^2$ be a closed oriented surface and let $\vec{\Phi}$ be a smooth immersion of $\Sigma^2$ into $\mathbb{R}^3$. Then $\vec{\Phi}$ is called \emph{Willmore immersion} if it is critical point of the Willmore functional $W$, i.e. if
\begin{displaymath}
\forall \vec{\xi} \in C^{\infty}_0 (\Sigma^2, \mathbb{R}^3) \; \; \;\frac{d}{dt} W(\vec{\Phi} + t \vec{\xi})_{t=0}=0.
\end{displaymath}
\end{Definition}

The corresponding Euler-Lagrange equation is due to Weiner \cite{Wei} and takes the form
\begin{equation}\label{eq:Weiner}
\Delta_g H +2H(H^2 -K)=0
\end{equation}
where $K$ is the Gauss curvature.

Unfortunately, there is a huge drawback in this equation: it is supercritical, i.e. while the functional can be defined for a  weak immersion (see later on for more details) with $L^2$-bounded second fundamental form, the corresponding Euler-Lagrange equation written as in \eqref{eq:Weiner} needs that $H$ is in $L^3$ to make sense, even in a distributional way.  But thanks to the significant break through in \cite{Riv2}, the equation can be reformulated as a conservation law that makes sense for weak immersions with $L^2$-bounded second fundamental form, i.e.  the problem becomes critical. Namely it holds (see \cite{Riv2})
\begin{Corollary}
A conformal immersion $\vec{\Phi}$ 
is Willmore if and only if 
\begin{equation}\label{EqWillmore}
div ( 2\nabla \vec{H} -3 H \nabla \vec{n} - \nabla^{\perp} \vec{n} \times \vec{H} ) =0
\end{equation}
where $\nabla^{\perp}f = (-\partial_{x_2}f, \partial_{x_1}f)$ denotes the rotated gradient.
\end{Corollary}

The fact that in the above corollary we assume the immersion to be conformal is not restrictive as we will see later on. \\

Before passing to the description of the  variational framework in which we will work,  let us recall the 
following lower bound on the Willmore energy due to Li and Yau  \cite{LY}.

\begin{theorem}\label{LiYau}
Let $\Sigma^2$ be a closed surface and let $\vec{\Phi}$ be a smooth immersion of $\Sigma^2$ into $\mathbb{R}^3$. Assume that there exists a point $p\in \mathbb{R}^3$ with at least $k$ pre-images under $\vec{\Phi}$. Then the following estimate holds
\begin{displaymath}
W(\vec{\Phi}) \geq 4\pi k.
\end{displaymath}
\end{theorem}

This estimate has a generalization to the context of varifolds which reads as follows (see e.g. 
Appendix of the article of Kuwert and Sch\"atzle \cite{KS3})
\begin{equation}\label{WillmoreInf}
\inf_{\Sigma^2 \; \textrm{smooth}}W(\Sigma^2)=4 \pi = \inf_{\mu \neq 0}W(\mu).
\end{equation}

The above Theorem \ref{LiYau} has the following useful corollary.
\begin{Corollary}\label{cor:embedding}
Let $\Sigma^2$ be a closed surface and let $\vec{\Phi}$ be a smooth immersion of $\Sigma^2$ into $\mathbb{R}^3$. If
\begin{displaymath}
W(\vec{\Phi}) < 8 \pi
\end{displaymath}
then $\vec{\Phi}$ is an embedding.
\end{Corollary}

\subsection{The variational framework}

The variational framework in which we will study the minimization problem is the one  of Lipschitz 
immersions with $L^2$-bounded second fundamental form that we are now going to discuss. \\

We start by briefly recalling the definition of the Sobolev spaces  $W^{k,p}(\Sigma^2,\mathbb{R}^3)$.\\ 
Let $\Sigma^2$ be a smooth closed oriented 2-dimensional manifold and let $g_0$ be a smooth reference metric on it.\\
Then the Sobolev spaces $W^{k,p}(\Sigma^2,\mathbb{R}^3)$ are defined as 
\begin{displaymath}
W^{k,p}(\Sigma^2,\mathbb{R}^3):= \Big\{f:\Sigma^2 \rightarrow \mathbb{R}^3 \; \Big| \; \sum_{j=1}^k \int_{\Sigma^2} \vert \nabla^j f \vert_{g_0}^p \; dvol_{g_0}< \infty\Big\}.
\end{displaymath}
Note that due to the fact that $\Sigma^2$ is assumed to be compact, the above spaces are independent of the metric $g_0$.\\
One particular case is of great interest for us: the  Lipschitz immersions $\vec{\Phi} \in W^{1,\infty} 
(\Sigma^2,\mathbb{R}^3)$ such that
\begin{equation}\label{nondeg}
\exists \; C \geq  1 \; \textrm{such that} \; C^{-1} g_0(X,X) \leq  \big|d\vec{\Phi} \cdot X \big|^2_{g_{\R^3}} \leq C g_0(X,X). 
\end{equation}

Next, observe that for a Lipschitz immersion satisfying \eqref{nondeg},  we have the following expression for the normal vector (called also Gauss map)
\begin{displaymath}
\vec{n}_{\vec{\Phi}}:=:=\frac{\partial_{x^i} \vec{\Phi} \times \partial_{x^j} \vec{\Phi} }{|\partial_{x^i} \vec{\Phi} \times \partial_{x^j} \vec{\Phi}|}\quad , 
\end{displaymath}
where $\vec{v_1}\times \vec{v_2}$ is the usual  cross product of vectors in $\R^3$. Notice that, from the definition, $\vec{n}$ is an $L^{\infty}$-map on $\Sigma^2$ with values in $\R^3$.\\

The space in which we will study our minimization problem is given by the Lipschitz immersions of $\Sigma^2$ with Gauss map in $W^{1,2}$ (or, equivalently, with $L^2$-bounded second fundamental form) defined as 
\begin{displaymath}
\mathcal{E}_{\Sigma^2}:= \left \{
\begin{split}
& \vec{\Phi} \in W^{1,\infty}(\Sigma^2,\mathbb{R}^3 )\; \textrm{such that} \\
& (\ref{nondeg}) \; \textrm{holds for some $C\geq 1$ and} \;
\int_{\Sigma^2} \vert d\vec{n}\vert^2_{g} dvol_{g} < \infty
\end{split}
\right \} \quad .
\end{displaymath}

Observe that the Willmore functional is well defined on the 
space $\mathcal{E}_{\Sigma^2}$, see  \eqref{WillmoreFunctional}. Moreover, thanks to the (local) existence of a smooth conformal structure (see Theorem \ref{conf} below), we can extend the notion of \emph{Willmore immersion} to the larger class $\mathcal{E}_{\Sigma^2}$ as follows (in the case of immersions into $\mathbb{R}^3$).

\begin{Definition}[Weak Willmore immersion]
Let $\vec{\Phi}$ belong to $\mathcal{E}_{\Sigma^2}$.\\
$\vec{\Phi}$ is called a \emph{weak Willmore immersion} if in any Lipschitz
conformal chart $\Psi$ from the two-dimensional disk $D^2$ into 
$(\Sigma^2,\vec{\Phi}^*g_{\mathbb{R}^3})$ the following holds
\begin{displaymath}
div ( 2\nabla \vec{H} -3 H \nabla \vec{n} - \nabla^{\perp} \vec{n} \times \vec{H} ) =0 \;
\textrm{in} \; \mathcal{D}'(D^2) .
\end{displaymath}
\end{Definition}

\subsection{The isoperimetric constraint}

As already mentioned in the introduction we want to study sequences of immersions of a closed surface $\Sigma^2$ of any genus minimizing the Willmore energy 
under the constraint that the isoperimetric ratio is fixed. More precisely, we minimize the Willmore energy among the weak immersions $\vec{\Phi}\in \mathcal{E}_{\Sigma^2} $ satisfying
\begin{displaymath} \label{isop}
iso(\vec{\Phi}):=\frac{(Area(\vec{\Phi}))^3}{\Big (Vol\Big( \; \overline{\vec{\Phi}} \;\Big)\Big)^2}=R \quad,
\end{displaymath}
where $R$ is a given constant and $Vol\Big( \; \overline{\vec{\Phi}} \; \Big)$ denotes the volume enclosed by $\vec{\Phi}$, i.e. the volume of the bounded
connected component of $\mathbb{R}^3 \backslash \vec{\Phi}(\Sigma^2)$. \\
Since by definition of $I_{\g}$ as in \eqref{eq:defI} we work below the energy threshold of $8 \pi$, thanks to Corollary \ref{cor:embedding} (or more precisely to its generalization to integer varifolds by Kuwert-Sch\"atzle in \cite[Appendix]{KS3}), the weak immersions we will handle do not have any self intersections. Thus there are no difficulties in defining the isoperimetric ratio.  \\

Observe that due to the scaling invariance of the Willmore functional,  fixing the  isoperimetric ratio is equivalent to fix 
both  area and enclosed volume.

At this point, it may be of interest to ask which equation is satisfied by a critical point of the Willmore functional under the  isoperimetric  constraint. By scaling, we can assume that the enclosed volume equals to one, so fixing the isoperimetric ratio is equivalent to fix the area; therefore we are looking for critical 
points of the Willmore functional under area constraint. Applying the Lagrange multiplier principle, we find the 
following characterization - for conformal immersions  which are area-constraint Willmore-
$$ 
\Delta_g H\vec{n}+2H(H^2 - K)\vec{n} =\mu H \vec{n} = \mu \vec{H} = \frac{\mu}{2} \Delta_g \vec{\Phi} = \frac{\mu}{2}e^{-2\lambda} \Delta \vec{\Phi} \quad,
$$
where we used $\vec{H}=\frac{1}{2} \Delta_g \vec{\Phi}=\frac{1}{2}e^{-2\lambda}\Delta \vec{\Phi}$.
Equivalently, we have
\begin{equation}\label{WillmoreConstrained}
 div ( 2\nabla \vec{H} -3 H \nabla \vec{n} - \nabla^{\perp} \vec{n} \times \vec{H}- \mu \nabla \vec{\Phi} )=0
\end{equation}
where we used the description of Willmore surfaces given by equation (\ref{EqWillmore}).\\
More details about this equation will be given later  when studying the regularity of the minimizer.

\subsection{Coulomb frames, conformal coordinates and estimates for the conformal factor}

As mentioned above, in addition to the formulation of the Willmore surface equation as a conservation law, we have to break the symmetry group, 
i.e. the invariance under diffeomorphisms of the surface $\Sigma^2$.\\
In fact, not only the right choice of coordinates, or gauge, is crucial, but in addition we need an appropriate control of the conformal factor.\\
A detailed presentation of the relation between conformal coordinates, Coulomb gauges and the Chern moving frame method is obviously beyond the scope of this article. Nevertheless, let us
shortly summarize the most important results.

First of all, the existence of conformal coordinates and a smooth conformal structure is asserted by the following theorem (see \cite{Riv4}, \cite{Hel})

\begin{theorem}\label{conf}
Let $\Sigma^2$ be a closed smooth 2-dimensional manifold. Let $\vec{\Phi}$ be an element of $\mathcal{E}_{\Sigma^2}$. Then there exists a finite covering of
$\Sigma^2$by discs $(U_i)_{i\in I}$ and Lipschitz diffeomorphisms $\Psi_i$ from the 2-dimensional unit disc $D^2$ into $U_i$ such that 
$\vec{\Phi} \circ \Psi_i$ realizes a Lipschitz conformal immersion of $D^2$. Since $\Psi_j^{-1} \circ \Psi_i$ is conformal and positive (i.e. holomorphic) on $\Psi_i^{-1}(U_i \cap U_j)$
the system of charts $(U_i, \Psi_i)$ defines
a smooth conformal structure $c$ on $\Sigma^2$ and in particular there exists a constant scalar curvature metric $g_c$ on $\Sigma^2$
and a Lipschitz diffeomorphism $\Psi$ of $\Sigma^2$ such that $\vec{\Phi} \circ \Psi$ realizes a conformal immersion of the Riemann surface 
$(\Sigma^2, g_c)$. 
\end{theorem}

Once we know about the existence of a conformal immersion, we want to address the question whether the conformal factor can be estimated appropriately. 
Concerning this latter question, we have the following answer (see \cite{RivNotes} and \cite{Riv5}).

\begin{theorem}
Let $\vec{\Phi}$ be a Lipschitz conformal immersion from the disc $D^2$ into $\mathbb{R}^3$. Assume 
\begin{equation}\label{eq:EnergyBound}
\int_{D^2} \vert \nabla \vec{n}_{\vec{\Phi}} \vert^2 < 8 \pi /3.
\end{equation}
Then for any $0< \rho <1$ there exists a constant $C_{\rho}$ independent of $\vec{\Phi}$ such that for the conformal factor $e^{\lambda}$ 
the following estimate holds
\begin{displaymath}
\sup_{p \in B^2_{\rho}(0)} e^{\lambda}(p) \leq C_{\rho} \Big \lbrack Area( \vec{\Phi}(D^2)) \Big \rbrack^{1/2}
\exp \Big ( C \int_{D^2} \vert \nabla \vec{n}_{\vec{\Phi}} \vert^2  \Big ).
\end{displaymath}
Moreover, for two given distinct points $p_1$ and $p_2$ in the interior of $D^2$ and again for $0 < \rho < 1$ there exists a constant $C_{\rho} > 0$
independent of $\vec{\Phi}$ such that
\begin{eqnarray} 
\vert \vert \lambda \vert \vert_{L^{\infty}(B^2_{\rho}(0))} &\leq& C_{\rho} \int_{D^2} \vert \nabla \vec{n}_{\vec{\Phi}} \vert^2  \nonumber \\ 
 &+&\, C_{\rho} \Big \vert \log \; \frac{\vert \vec{\Phi}(p_1) - \vec{\Phi}(p_2) \vert }{\vert p_2 - p_1 \vert} \Big \vert 
\,+ \, C_{\rho} \, \log^+ \big \lbrack C_{\rho} Area( \vec{\Phi}(D^2)) \big  \rbrack \nonumber
\end{eqnarray}
where $\log^+$ denotes the positive part of the logarithm. 
\end{theorem}

In our particular case, we will have a variant of such an estimate of the conformal factor. \\
For a precise statement of this modified estimate, we refer to Lemma \ref{confFactor} below. \\

Note that the (Willmore) energy assumption \eqref{eq:EnergyBound} is also crucial in H\'elein's  lifting theorem which asserts the existence of a moving frame with energy estimates provided that one starts with a map whose second fundamental form, measured in $L^2$, is below the critical threshold.

\begin{theorem}\label{Helein}
Let $\vec{n}$ be a $W^{1,2}$-map from the 2-dimensional disc $D^2$ into the unit sphere ${\mathbb S}^2\subset \R^3$. Then there exists a constant  $C>0$ such that there exist $\vec{e}_1$ and $\vec{e}_2$ in $W^{1,2}(D^2,{\mathbb S}^2)$ such that
\begin{displaymath}
\vec{n}=*(\vec{e}_1 \wedge \vec{e}_2)
\end{displaymath}
with
\begin{displaymath}
\int_{D^2} \sum_{i=1}^{2}\vert \nabla \vec{e}_i\vert^2 < C\int_{D^2} \vert \nabla \vec{n}\vert^2 
\end{displaymath}
provided that
\begin{displaymath}
\int_{D^2} \vert \nabla \vec{n}\vert^2 < \frac{8 \pi}{3}.
\end{displaymath}
\end{theorem}
A proof of this theorem can be found in \cite{Hel}.

\section{Proof of Theorem \ref{thm:Schy}}
First of all, let us state the following  useful result proved by Schygulla (see\cite{Schy}).

\begin{theorem}\label{Pro:Schy}
For every $R\in[36\pi, +\infty)$ there exists a smooth embedded spherical surface  with isoperimetric ratio equal to $R$ and having Willmore energy strictly less than $8\pi$.
\end{theorem} 

Let us recall the clever arguments in \cite{Schy} to prove the theorem above: Inverting a Cathenoid in the origin and desingularizing it in $0$, one obtains a spherical surface with energy strictly less than $8 \pi$ and arbitrarily small isoperimetric ratio. Then, let this surface evolve under Willmore flow. Thanks to the results of Kuwert-Sch\"atzle (see \cite{KS1}, \cite{KS2} and \cite{KS3}), the flow will converge smoothly (as $t\uparrow +\infty$) to a round sphere (whose isoperimetric ratio is $36\pi$). Therefore, since the flow is smooth and does not increase the Willmore energy and since the isoperimetric ratio depends continuously on the parameter of the flow, for every $R\in [36\pi, +\infty)$ one has produced a spherical surface with isoperimetric ratio equal to $R$ and Willmore energy strictly less than $8 \pi$ as desired.
\\

The second ingredient for the proof of Theorem \ref{thm:Schy} is a compactness result for weak immersions of spheres proved by the second and third authors \cite[Theorem I.2]{MR}. Since by Theorem \ref{Pro:Schy} we can work under an $8 \pi-\delta$ assumption, let us state the compactness result under this simplifying hypothesis (which prevents the bubbling phenomenon and the presence of branch points, thanks to Theorem \ref{LiYau}).

\begin{theorem}\label{Pro:MR}
Let $\{\vec{\Phi}_k\}_{k \in \N}\subset {\mathcal E}_{\mathbb{S}^2}$ be a sequence of  weak  conformal immersions  of the $2$-sphere $\mathbb{S}^2$ into $\R^3$ such that
\begin{equation}
\label{I.4}
\sup_{k} Area(\vec{\Phi}_k) <\infty,\quad \inf_k diam(\vec{\Phi}_k(\mathbb{S}^2))>0 \; \text{ and } W(\vec{\Phi}_k)<8\pi-\delta \quad,
\end{equation}
 where $diam(E)$ is the diameter of the subset $E\subset \R^3$ and $\delta>0$ is some positive constant.\\
  
Then there exist a subsequence that we still denote $\{\vec{\Phi}_k\}_{k \in \N}$,  a sequence $\{f_k\}_{k \in \N}$ of elements in ${\mathcal M}^+(\mathbb{S}^2)$ (the positive M\"obius group of $\mathbb{S}^2$) and  finitely many points $\{a_1, \ldots, a_n\}$ such that
\begin{equation} \label{I.6}
\vec{\Phi}_k\circ f_k\rightharpoonup \vec{\Phi}_\infty\quad\text{ weakly in } W^{2,2}_{loc}(\mathbb{S}^2\setminus\{a_1\cdots a_n\})\quad,
\end{equation}
where $\vec{\Phi}_\infty \in{\mathcal E}_{\mathbb{S}^2}$ is conformal.  In addition, by lower semicontinuity of  $W$ under weak $W^{2,2}$ convergence, we have
\begin{equation}\nonumber
W(\vec{\Phi}_\infty) \leq  \liminf_{k} W(\vec{\Phi}_k)< 8\pi \quad,
\end{equation}
so $\vec{\Phi}_\infty$ is a (weak) embedding, and moreover 
\begin{equation}
\label{I.6b}
Area(\vec{\Phi}_k)\to  Area (\vec{\Phi}_\infty)\quad\text{and} \quad Vol\Big( \; \overline{\vec{\Phi}_k} \; \Big) \to Vol\Big( \; \overline{\vec{\Phi}_\infty} \; \Big). 
\end{equation}
\end{theorem}

The proof Theorem \ref{thm:Schy} now follows quite easily by the two theorems above. Thanks to Theorem \ref{Pro:Schy}, for any $R\in[36\pi,+\infty)$, the infimum of $W$ among weak immersions in ${\mathcal E}_{\mathbb{S}^2}$ under the constraint of fixed isoperimetric ration equal to $R$ is strictly less than then $8 \pi$. Therefore, for any minimizing sequence $\{\vec{\Phi}_k\}_{k \in \N}\subset {\mathcal E}_{\mathbb{S}^2}$ of the constrained problem we have
$$W(\vec{\Phi}_k)\leq 8\pi-\delta \quad, $$
for some $\delta>0$. Therefore, the $\vec{\Phi}_k$ are (weak) embeddings and, as  already observed by the scale invariance of $W$, we can assume that the enclosed volume of $\vec{\Phi}_k$ is constantly equal to $1$ and the area of $\vec{\Phi}_k$ is constantly equal to $R^{1/3}$. From Simon's Lemma \ref{LemmaSimon} recalled in the Appendix, we also have a strictly positive lower bound on the diameters of $\vec{\Phi}_k(\mathbb{S}^2)$ as subsets of $\R^3$. Collecting the above informations, we conclude that the $\vec{\Phi}_k$ satisfy the assumptions of Theorem \ref{Pro:MR}. It follows that there exists a weak embedding $\vec{\Phi}_\infty\in  {\mathcal E}_{\mathbb{S}^2}$, with $iso(\vec{\Phi}_\infty)=R$, which realizes the infimum of $W$ under the isoperimetric constraint, as desired. For the proof of the regularity see Subsection \ref{Subsub:reg}. \hfill$\Box$
\\

In the next section we will analyze the more delicate case of genus $\g\geq 1$.

\section{Proof of Theorem \ref{MainTheorem}}

Let $\Sigma^2$ be a  2-dimensional surface of genus $\g\geq 1$, and $\left \{ \vec{\Phi}_k \right \}_{k \in \mathbb{N}}\subset \mathcal{E}_{\Sigma^2}$ be a minimizing  sequence for the minimization problem of Theorem \ref{MainTheorem}, i.e.  minimization of the Willmore energy among weak immersions of $\Sigma^2$ under the constraint that
the isoperimetric ratio $\frac{(Area(\vec{\Phi}_n))^3}{\Big(Vol\Big( \; \overline{\vec{\Phi}_n} \;\Big)\Big)^2}$ is fixed and equals $R$, for some $R \in I_{\g}$ where $I_{\g}$ is defined in \eqref{eq:defI}. Our goal is to show that the minimizing sequence is compact, the  limit belongs to $\mathcal{E}_{\Sigma^2}$ and satisfies the required geometrical properties.  The procedure we apply here is inspired by the one used in \cite{Riv5} (see also \cite{BeRi2}, Lemma II.1).

\subsection{Normalization of the minimizing sequence}

\subsubsection{Conformality of $\vec{\Phi}_k \in \mathcal{E}_{\Sigma^2}$}

Note that a priori, the elements of our sequence are not necessarily conformal. But thanks to theorem \ref{conf} we may replace the original $\vec{\Phi}_k$s by $\vec{\Phi}_k \circ \Psi_k$ which are weakly conformal immersions.\\
In a few words, the idea behind the quoted theorem is that one starts with a frame with controlled energy, improves it into a Coulomb frame on each coordinate patch, applies the improvement of the well known Wente estimate due to Chanillo and Li (see \cite{CL}) and finally concludes by using the Riemann mapping theorem.\\
This new sequence will still be denoted - by abuse of the notation - $\left \{ \vec{\Phi}_k \right \}_{k \in \mathbb{N}}$. Note that this new sequence consists of conformal immersions of $\Sigma^2$. 
Rephrased, we have exploited the \emph{invariance in the domain}. \\
Note that this procedure does not affect the isoperimetric ratio since of course this latter quantity is intrinsic, i.e. it does not depend on the choice of coordinates. \\

{\bf Remark}
At this stage, we would briefly comment on the difference between the present procedure and the one used in the unconstrained case presented in \cite{Riv5}. \\
In this latter case the \textit{invariance in the target under M\"obis transformations} is used as well (see \cite{Riv5}, Lemma A.4). In particular, inversions play a crucial role. Since such inversions do not preserve the isoperimetric ratio we can and will not exploit this invariance.

\subsubsection{Points of energy concentration}

Roughly speaking, we will have only finitely many points $\left \{ a_i \right \}$ where energy can accumulate with a critical energy threshold  equal to $\frac{8 \pi}{3}$. This is done as follows. \\
First of all, recall that $g_k:= \vec{\Phi}_k^* g_{\mathbb{R}^3}$ is the pull back metric, and $h_k$ denotes the metric of constant scalar curvature  associated to the conformal class $c_k$ of the metric $g_k$. \\
Next,  the hypothesis that the Willmore energy of our immersions stays below $\min\{8 \pi,\omega^3_{\g}\}$ (for sufficiently large $k$) implies that the conformal structures are contained in a compact subset of the moduli space of $\Sigma^2$ (see Theorem \ref{ConfClass} in the Appendix). Thus, the metrics $h_k$ converge to $h_{\infty}$, the metric of constant scalar curvature associated to the limiting conformal structure $c_\infty$ (about the existence of such a metric we refer to \cite{J}). This convergence holds in $C^s(\Sigma^2)$ for all $s$. \\ 
Now, to each point $x$ in $\Sigma^2$ we assign a critical radius which ``cuts out $\frac{8 \pi}{3}$ Willmore energy'', more precisely, this radius $\rho^k_x$ is defined as follows
\begin{displaymath}
\rho^x_k:= \inf \Big\{ \rho \; \Big\vert \; \int_{B_{\rho_k}(x)} \vert \nabla n_{\vec{\Phi}_k} \vert_{g_k}^2 dvol_{g_k}=
\int_{B_{\rho_k}(x)} \vert \nabla n_{\vec{\Phi}_k} \vert_{h_k}^2 dvol_{h_k}=  
\frac{8 \pi}{3} \Big\}
\end{displaymath}
where $B_{\rho_k}(x)$ denotes the geodesic ball in $(\Sigma^2, g_k)$ centered at $x$ with radius $\rho_k$.
Obviously, $\{ B_{\rho^x_k/2}(x) \}_{x \in \Sigma^2}$ is a covering of  $\Sigma^2$. \\
Next, we extract a finite Besicovitch covering: each point $x \in \Sigma^2$ is covered by at most $N=N(\Sigma^2, g_{\infty}) \in \mathbb{N}$ such balls. This extracted covering is denoted by $\{ B_{\rho^i_{k}/2}(x^i_k) \}_{i \in I}$. 
Then we pass to a subsequence such that the following properties are satisfied
\begin{itemize}
\item[i)]
\begin{displaymath}
I  \quad \textrm{is independent on} \; k
\end{displaymath}
\item[ii)]
\begin{displaymath}
x_k^i \rightarrow x_{\infty}^i
\end{displaymath}
\item[iii)]
\begin{displaymath}
\rho_{k}^i \rightarrow \rho_{\infty}^i.
\end{displaymath}
\end{itemize}
And we set 
\begin{displaymath}
I_0:= \left \{ i \in I \; \textrm{such that} \; \rho_{\infty}^i =0 \right \} \; \; \textrm{and} \; \; I_1 := I \backslash I_0.
\end{displaymath}
Obviously, 
\begin{displaymath}
\bigcup_{i \in I_1}\bar{B}_{\rho_{\infty}^i/2}(x^i_{\infty}) \; \textrm{covers} \; \Sigma^2
\end{displaymath}
(where the balls are measured in the metric $g_{\infty}$)
and thanks to the strict convexity of our balls - with respect to the  flat or hyperbolic metric - there are only isolated, finitely many points in $\Sigma^2$ which are not covered by the union of the open balls. Thus, we denote these exceptional point by $\{ a_i, \dots , a_n \}$, or more precisely
\begin{displaymath}
\{ a_1, \dots, a_n \} := \Sigma^2 \backslash \bigcup_{i \in I_1} B_{\rho_{\infty}^i/2}(x^i_{\infty}).
\end{displaymath} 
Rephrased, we have identified the points of energy concentration.\\

Note that this last covering of $\Sigma^2 \backslash \{ a_1, \dots, a_n \}$,
\begin{equation}\label{cover}
\Sigma^2 \backslash \{ a_1, \dots, a_n \} \subset \bigcup_{i \in I_1}B_{\rho_{\infty}^i/2}(x^i_{\infty})
\end{equation}
 satisfies
\begin{displaymath}
\int_{B_{\rho_{\infty}^i}(x^i_{\infty})} \vert \nabla \vec{n}_{\vec{\Phi}_k} \vert^2 \leq \frac{8 \pi}{3}
\end{displaymath}
for all $k$ and all $i \in I_1$.

\subsubsection{Analysis away from the energy-concentration-points: Control of the conformal factor}

In this subsection we will perform the analysis of our sequence $\{ \vec{\Phi}k \}_{k \in \mathbb{N}}$ away from the energy-concentration points $\{a_1, \dots , a_n \}$. Let us start with an estimate for the conformal factors. As above, denote the conformal factor of $\vec{\Phi}_k$ by $\lambda_k$. 
\begin{Lemma}\label{confFactor}
Let $\varepsilon >0$. Then there exist constants $C_k$ and a constant $C_{\varepsilon}$ such that - up to passing to a subsequence -
we have
\begin{displaymath}
\vert \vert \lambda_k - C_k \vert \vert_{L^{\infty}(\Sigma^2 \backslash \cup^{n} _{i=1}B_{\varepsilon}(a_i))} \leq C_{\varepsilon}.
\end{displaymath}
\end{Lemma}
Note that $C_{\varepsilon}$ does not depend on $k$, but only on $\varepsilon$ whereas the $C_k$ may depend on $k$. \\

A result of the same spirit is used in \cite{Riv5}. Although the proof  is quite similar to the one in \cite{Riv5}, let us recall the main arguments for the reader's convenience. \\

\textit{Proof of lemma \ref{confFactor}:}\\
First of all, since $\Sigma^2$ is a surface of genus at least one, the reference metric $g_0$ of constant scalar curvature is  flat or hyperbolic. Therefore, recalling that the immersions $\vec{\Phi}_k$ are conformal so that  $g_k:=\vec{\Phi}_k^*g_{\mathbb{R}^3}= e^{2\lambda_k}g_0$, the conformal factors $\lambda_k$ satisfy the Gauss-Liouville equation

\begin{displaymath}
-\Delta_{g_0}\lambda_k = K_{g_k}e^{2\lambda_k} + K_{g_0} \quad .
\end{displaymath} 
Now, since the second fundamental form is in $L^2$, observe that $K_{g_k}e^{2\lambda_k}$ belongs to $L^1$ - with respect to $g_0$. This fact implies by standard elliptic estimates that for a constant depending only on the surface $\Sigma^2$ and the metric $g_0$ we have
\begin{equation}
\nonumber
\begin{split}
& \vert \vert d \lambda_k \vert \vert _{L^{2,\infty}(\Sigma^2,g_0)}\leq C \vert \vert \Delta_{g_0} \lambda_k \vert \vert_{L^1(\Sigma^2,g_0)} 
\leq C \left[ \int_{\Sigma^2} \vert K_{g_k} \vert e^{2\lambda_k} dvol_{g_0} + |K_{g_0}| Area_{g_0}(\Sigma^2) \right] \\
& \quad \leq C' \left[  \int_{\Sigma^2} \vert K_{g_k} \vert dvol_{g_k} +1\right] \leq C' \left[ \int_{\Sigma^2} \vert d \vec{n}_{\vec{\Phi}_k} \vert_{g_k}^2 dvol_{g_k}+1 \right]\\
& \quad \leq C'' \left[W(\vec{\Phi}_k)+1\right]< C'''.
\end{split}
\end{equation} 
Now, let $\varepsilon$ be given. Starting from the covering $\bigcup_{i \in I_1}B_{\rho_{\infty}^i/2}(x^i_{\infty})$ (cf. (\ref{cover})) we 
obtain a covering of $\Sigma^2 \backslash \cup^{n} _{i=1}B_{\varepsilon}(a_i)$ of the following form 
\begin{displaymath}
\Sigma^2 \backslash \cup^{n} _{i=1}B_{\varepsilon}(a_i) \subset \bigcup_{i \in I_1}B_{r^i/2}(x^i_{\infty})
\end{displaymath}
where $r^i < \rho_{\infty}^i$ (the balls here and in the following steps are measured in the metric $g_{\infty}$). \\
Note that the connectedness of $\Sigma^2$ allows us, up to a relabeling of our balls, to assume that two consecutive balls have non-empty intersection.
Recall now that by the properties of the covering (\ref{cover}) and the conformal invariance of the integrand, we have
\begin{displaymath}
\int_{B_{r^i}(x_{\infty}^i)}\vert \nabla \vec{n}_{\vec{\Phi}_k}\vert^2 < \frac{8 \pi}{3} \quad .
\end{displaymath}
Next ,we  use H\'eleins moving frame method \cite{Hel}: \\
Thanks to our construction, on the balls $B_{r^i}(x^i_{\infty})\subset B_{\rho_{\infty}^i}(x^i_{\infty})$ we have strictly less than the critical  energy $\frac{8 \pi}{3}$. Using Theorem \ref{Helein} (for more details see \cite{Hel} or \cite{RivNotes}), upon identifying $B_{r_i}(x_{\infty}^i)$ with the 2-dimensional Euclidean unit disk $D^2$, for each $k$ there exists a moving frame with controlled energy, i.e. $(\vec{e}^1_k,\vec{e}^2_k) \in {\mathbb S}^2 \times {\mathbb S}^2$ with the following properties
\begin{itemize}
\item[i)]
\begin{displaymath}
\vec{e}^{\,1}_k \cdot \vec{e}^{\,2}_k=0 \; \textrm{and} \; \vec{n}_{\vec{\Phi}_k}= \vec{e}^{\,1}_k \wedge \vec{e}^{\,2}_k
\end{displaymath}
\item[ii)]
\begin{displaymath}
\int_{D^2}( \vert \nabla \vec{e}^{\,1}_k\vert^2 +  \vert \nabla \vec{e}^{\,2}_k\vert^2 ) \leq 2 \int_{D^2}\vert \nabla \vec{n}_{\vec{\Phi}_k}\vert^2 < \frac{16 \pi}{3}
\end{displaymath}
\item[iii)]
\begin{displaymath}
\left \{ 
\begin{array}{ll}
div(\vec{e}^{\,1}_k, \nabla \vec{e}^{\,2}_k)  & =  0  \quad \textrm{in $D^2$} \\
\Big ( \vec{e}^{\,1}_k, \frac{\partial \vec{e}^{\,2}_k}{\partial \nu} \Big )  & = 0 \quad \textrm{on $\partial D^2$}.
\end{array} 
\right .
\end{displaymath}
\end{itemize}
In these frames we can express the conformal factors $\lambda_k$ as follows
\begin{displaymath}
-\Delta \lambda_k = (\nabla ^{\perp} \vec{e}^{\,1}_k, \nabla \vec{e}^{\,2}_k).
\end{displaymath}
Now, consider the solution $\mu_k$ of the following problem
\begin{displaymath}
\left \{ 
\begin{array}{rll}
- \Delta \mu_k & =  (\nabla ^{\perp} \vec{e}^{\,1}_k, \nabla \vec{e}^{\,2}_k)  &  \textrm{in $D^2$} \\
\mu_k  & = 0 & \textrm{on $\partial D^2$}.
\end{array} 
\right .
\end{displaymath}
For this solution, Wente's theorem (see theorem \ref{WenteEtAl} in the Appendix) gives the estimate
\begin{equation}
\nonumber
\begin{split}
& \vert \vert \mu_k \vert \vert_{L^{\infty}(B_{r^i}(x^i_{\infty}))} + \vert \vert \nabla \mu_k \vert \vert _{L^{2,1}(B_{r^i}(x^i_{\infty}))} + \vert \vert \nabla^2 \mu_k \vert \vert_{L^1(B_{r^i}(x^i_{\infty}))} \\
& \quad \leq C \int_{B_{r^i}(x^i_{\infty} )}( \vert \nabla \vec{e}^{\,1}_k\vert^{\,2} +  \vert \nabla \vec{e}^2_k\vert^2 ) 
\leq 2 C \int_{B_{r^i}(x^i_{\infty} )}\vert \vec{n}_{\vec{\Phi}_k}\vert^2 < C \frac{16 \pi}{3}. 
\end{split}
\end{equation}
Next, we look at $v_k:= \lambda_k - \mu_k$. Since $v_k$ is harmonic, then we have (see e.g. \cite{GT})
\begin{displaymath}
\vert \vert v_k -\bar{v}_k \vert \vert _{L^{\infty}(B_{r^i/2}(x^i_{\infty}))} \leq C
\end{displaymath}
where $\bar{v}_k$ is the average of $v_k$ over the ball $B_{r^i}(x^i_{\infty})$. 
Putting together all the information we have so far, we conclude  that there exist constants $\bar{C}_k$
\begin{displaymath}
\vert \vert \lambda_k - \bar{C}_k \vert \vert_{L^{\infty}(B_{r^i/2}(x^i_{\infty}))} \leq C.
\end{displaymath}
In a last step, we combine the fact that the above arguments apply for all balls in our finite covering with the fact that
two consecutive balls have non-empty intersection in order to conclude that the constant $\bar{C}_k$ depends only on $k$
but not on the ball we look at.
This completes the proof of the lemma. \hfill$\Box$
\\

Now, a priori there are three possibilities: either $C_k$ remain bounded, or they tend to $-\infty$ or  diverge to $+\infty$.\\
Observe that The latter case is excluded, since our hypothesis imply that the area of $\vec{\Phi}_k$ remain bounded. Indeed, since the Willmore functional is scaling invariant, the isoperimetric constraint is equivalent to fix both area and
enclosed volume. Thus, recalling that the area form is exactly
$e^{2\lambda_k}$, the assumption that $C_k \rightarrow +\infty$ would imply that the area of $\vec{\Phi}_k(\Sigma^2)$, which is of course 
at least the area of $\Sigma^2 \backslash \cup^{n} _{i=1}B_{\varepsilon}(a_i))$, would become arbitrarily large. But this contradicts the assumption that the areas are fixed. \\
On the other hand, the fact that if $C_k$ tend to $- \infty$ the area of $\Sigma^2 \backslash \cup^{n} _{i=1}B_{\varepsilon}(a_i))$
tends to zero does not a priori lead to a contradiction to the boundedness of the areas. In fact, due to a bubbling phenomenon, it could happen that we have large area on the 
exceptional balls $B_{\varepsilon}(a_i)$. So, a priori we can exclude only the possibility that $C_k \rightarrow + \infty$.\\
Below we will examine the remaining  two possibilities.\\

\subsection{The case of bounded conformal factors}\label{Subsec:CkBounded}
In the present  subsection, we assume that the $C_k$ remain bounded. The case of diverging conformal factors will be discussed (and excluded) in  Subsection \ref{Subsec:Diverging}.

\subsubsection{Weak convergence in $W^{2,2}$}

The assumption that the $C_k$ remain bounded, i.e. $\sup_{k \in \mathbb{N}} \vert C_k \vert < \infty$, together with Lemma \ref{confFactor} immediately implies that
\begin{displaymath}
\limsup_{k \rightarrow \infty} \vert \vert \lambda_k \vert \vert_{L^{\infty}(\Sigma^2\backslash \cup^{n} _{i=1}B_{\varepsilon}(a_i))} < \infty.
\end{displaymath}
Therefore, recalling that $\Delta \vec{\Phi}_k = 2 e^{2 \lambda_k}\vec{H}_k$, we infer that the following estimates hold
\begin{displaymath}
\left \{ 
\begin{array}{l}
\vert \vert \Delta \vec{\Phi}_k \vert \vert_{L^2(\Sigma^2\backslash \cup^{n} _{i=1}B_{\varepsilon}(a_i))} < C(\varepsilon)\\ 
\vert \vert \log \vert \nabla \vec{\Phi}_k \vert \, \vert \vert_{L^{\infty}(\Sigma^2\backslash \cup^{n} _{i=1}B_{\varepsilon}(a_i))} < C(\varepsilon)
\end{array} 
\right .
\end{displaymath}
From that, we deduce that there exists a subsequence which converges weakly to a limit $\vec{\Phi}_{\infty}$ in $W^{2,2}(\Sigma^2\backslash \cup^{n} _{i=1}B_{\varepsilon}(a_i))$; in particular, by Rellich Theorem, we have that  
\begin{displaymath}
\nabla \vec{\Phi}_k \rightarrow \nabla \vec{\Phi}_{\infty}  \;  \textrm{ strongly in} \; L^p(\Sigma^2 \backslash \cup_{i=1}^n B_{\varepsilon}(a_i)) \; \forall \, p < \infty
\end{displaymath}
and hence upon passing to a subsequence, the gradients converge even almost everywhere pointwise.\\

\begin{Remark}\label{rem:IsopLimit}
\rm{Note that the above strong convergence of the gradients implies the convergence of the areas - away from the points $a_i$.\\
Moreover, by Sobolev embeddings we know that the sequence $\vec{\Phi}_k$ converges in $C^{0,\gamma}$ for $\gamma < 1$ and thus, roughly said, the enclosed volume converges as well away the points $a_i$ and hence  ``the  isoperimetric ratio is preserved in the limit away the points $a_i$''. We will make this statement more precise later in Subsubsection \ref{Subsub:Isop}.}\hfill$\Box$
\end{Remark}

\subsubsection{Conformality of the limit}

Recall that the sequence we have after all the preceding steps consists of conformal embeddings, i.e. we have for all $k$
\begin{displaymath}
\left \{ 
\begin{array}{rl}
\partial_{x_1}\vec{\Phi}_k \cdot \partial_{x_2}\vec{\Phi}_k & =  0\\
\vert \partial_{x_1}\vec{\Phi}_k \vert^2 - \vert \partial_{x_2}\vec{\Phi}_k \vert^2 & = 0 .
\end{array} 
\right .
\end{displaymath} 

Due to the a.e. pointwise convergence of the gradients $\nabla \vec{\Phi}_k$ (see above)
we immediately can conclude that the conformality condition passes to the limit, i.e.
\begin{displaymath}
\left \{ 
\begin{array}{rl}
\partial_{x_1}\vec{\Phi}_{\infty}\cdot \partial_{x_2}\vec{\Phi}_{\infty} & =  0\\
\vert \partial_{x_1}\vec{\Phi}_{\infty} \vert^2 - \vert \partial_{x_2}\vec{\Phi}_{\infty} \vert^2 & = 0 .
\end{array} 
\right .
\end{displaymath} 

Together with the $L^{\infty}$-control of the conformal factors, this implies that 
\begin{displaymath}
\vec{\Phi}_{\infty} \in  \mathcal{E}_{\Sigma^2 \backslash \cup^{n} _{i=1}B_{\varepsilon}(a_i)}.
\end{displaymath}

\subsubsection{Control of $\vec{\Phi}_{\infty}$ over the whole $\Sigma^2$}

Thanks to the results of the analysis performed so far, we can apply the following lemma due to Rivi\`ere (see \cite{Riv5}, see also the article of Kuwert and Li  \cite{KL}).

\begin{Lemma}\label{LemmaLipschitz}
Let $\vec{\xi}$ be a conformal immersion of $D^2 \backslash \{ 0 \}$ into $\mathbb{R}^3$ in \\
$W^{2,2}(D^2\backslash \{ 0 \}, \mathbb{R}^3)$ and such that $\log  \vert \nabla \vec{\xi} \vert  \in L^{\infty}_{loc}(D^2 \backslash \{0 \})$. Assume $\vec{\xi}$ extends to a map in $W^{1,2}(D^2)$ and that the corresponding Gauss map $\vec{n}_{\vec{\xi}}$ also extends to a map in $W^{1,2}(D^2, {\mathbb S}^2)$. Then $\vec{\xi}$ realizes a Lipschitz conformal immersion of the whole disc $D^2$ and there exits a positive integer $n$ and a constant $C$ such that
\begin{displaymath}
(C- o(1)) \vert z \vert ^{n-1} \leq \Big \vert \frac{\partial \vec{\xi}}{\partial z} \Big \vert \leq 
(C+o(1)) \vert z \vert ^{n-1}.
\end{displaymath}
\end{Lemma}

More precisely, we apply this lemma to $\vec{\Phi}_{\infty}$ around each exceptional point $a_i$, $i=1, \dots, n$. \\

\textit{Claim: In our situation, the assertion of the above lemma holds with $n=1$, i.e.
\begin{displaymath}
(C- o(1)) \vert z \vert \leq \Big \vert \frac{\partial \vec{\xi}}{\partial z} \Big \vert \leq 
(C+o(1)) \vert z \vert.
\end{displaymath}
}

In other words, there is no branching.\\

\textit{Proof of the claim:}\\
First of all, observe that the above chain of inequalities implies that
for any $\delta > 0$ there exists a radius $r_{\delta}> 0$ such that for all $r< r_{\delta}$ we have
\begin{displaymath}
\vec{\Phi}_{\infty} (B_r(a_i)) \subset B_{\rho}(\vec{\Phi}_{\infty}(a_i)) \quad \textrm{and} \quad 
\vert \partial_z \vec{\Phi}_{\infty} \vert = e^{\lambda_{\infty}} \geq C \frac{1-\delta}{\sqrt{2}} \vert z \vert^{n-1}.
\end{displaymath}
where $\rho=C2^{-1/2}n^{-1}(1+ \delta)r^n$. \\
Using these facts, we then can estimate the mass of $\vec{\Phi}_{\infty}(\Sigma^2)$ inside the ball $B_{\rho}(\vec{\Phi}_{\infty}(a_i))$ 
as follows (note that here $\vec{\Phi}_{\infty}$ is seen as a varifold)
\begin{eqnarray}
\mu(\vec{\Phi}_{\infty} \cap B_{\rho}(\vec{\Phi}_{\infty}(a_i)) &\geq& C^2 \frac{(1- \delta)^2}{2} \int_{B_r(a_i)} \vert z \vert^{2n-2} \nonumber \\
&\geq& \frac{\pi C^2 (1- \delta)^2}{2n}r^{2n} \nonumber \\
&\geq& n \pi \Big (\frac{1- \delta}{1 + \delta} \Big)^2 \rho^2. \nonumber
\end{eqnarray}
From this estimate we deduce that the 2-dimensional lower density of $\vec{\Phi}_{\infty}(\Sigma^2)$ at the point $\vec{\Phi}_{\infty}(a_i)$,
$\theta^2_*(\vec{\Phi}_{\infty}(\Sigma^2),\vec{\Phi}_{\infty}(a_i))$ is bigger or equal to $n$.\\
On the other hand, the Li-Yau inequality (see \cite{LY}) - and in particular the extension of this inequality to the setting of varifolds 
with mean curvature in $L^2$ (see \cite{KS3} (Appendix))
\begin{displaymath}
\theta^2_*(\mu,x_*) \leq \frac{1}{4 \pi}W(\mu)
\end{displaymath}
implies that
\begin{displaymath}
n \leq \theta^2_*(\vec{\Phi}_{\infty}(\Sigma^2),\vec{\Phi}_{\infty}(a_i)) \leq \frac{W(\vec{\Phi}_{\infty}(\Sigma^2))}{4 \pi} 
< \frac{8 \pi - \delta}{4 \pi} < 2
\end{displaymath}
where we used also the lower semicontinuity 
of the Willmore functional $W$ and the assumption that $W(\vec{\Phi}_k) < 8 \pi- \delta$ 
for some $\delta > 0$.\\
Finally, this leads to the conclusion that $n=1$.
\hfill$\Box$ 
\\

In other words, we have a point removability phenomenon (cf. also proof of the above cited lemma). \\

Note that this together with the fact that away from the points $a_i$ we have convergence in 
$C^{0,\gamma}$ implies that $\vec{\Phi}_k(\Sigma^2)$ remain bounded.

\subsubsection{Limit of $\left \{ \vec{\Phi}_k \right \}_{k \in \mathbb{N}}$ in $\mathcal{E}_{\Sigma^2}$}

The results from the previous section immediately lead to the conclusion that $\vec{\Phi}_{\infty}$ is a Lipschitz immersion, and thus together
with the $W^{2,2}$-convergence established earlier we find that $\vec{\Phi}_{\infty} \in \mathcal{E}_{\Sigma^2}$.

\subsubsection{The isoperimetric constraint in the limit}\label{Subsub:Isop}
It remains to show that the limit $\vec{\Phi}_{\infty}$ satisfies the isoperimetric constraint. \\
First of all, observe that in all the modifications of our initial minimizing sequence, we did nothing that could affect the isoperimetric constraint.
So, we always have a sequence of immersions respecting the isoperimetric constraint.\\
Recall that the requirement of fixed isoperimetric ratio 
can be rephrased as fixing the area as well as the enclosed volume.\\
At first glimpse, one might have the idea that the information we have at hand - more precisely that 
our sequence $\{\vec{\Phi}_k\}_{k \in \N}$ converges in $W^{2,2}$ - away from points of energy concentration - 
might be enough in order to show that the limit 
$\vec{\Phi}_{\infty}$ satisfies the required isoperimetric constraint. But this is not the case.
We will explain this by looking at the area. \\
Convergence of the area would be a consequence of 
\begin{displaymath}
e^{\lambda_k} \rightarrow e^{\lambda_{\infty}}.
\end{displaymath}
But this convergence does not need to hold, since we have only the local control of the conformal factors
\begin{displaymath}
 \vert \vert \lambda_k - C_k \vert \vert_{L^{\infty}(\Sigma^2 \backslash \cup^{n} _{i=1}B_{\varepsilon}(a_i))} \leq C_{\varepsilon}
\end{displaymath}
from Lemma \ref{confFactor}. Thus, the closer we get to the points $a_i$ the bigger the $L^{\infty}$-norms
of the conformal factors may become. \\

Our strategy in order to show that the limit $\vec{\Phi}_{\infty}$ satisfies the isoperimetric constraint
is to exclude a bubbling phenomenon.\\ 
Roughly speaking, we will detect regions of positive area which carry an energy contribution
of at least $4 \pi$. This will lead to a contradiction to our initial hypothesis that 
$W(\vec{\Phi}_k)< 8\pi -\delta$. More precisely, we prove the following lemma.

\begin{Lemma}\label{lemmaBubbling}
Let $\{ \vec{\Phi}k \}_{k \in \mathbb{N}}$ be as above and denote by $a_i$, $i=1, \dots n$, the points of energy concentration.
Then we have the following assertions
\begin{itemize}
\item[i)]
Assume that there exists an index $i$ such that
\begin{displaymath}
\liminf_{\varepsilon \rightarrow 0} \liminf_{k \rightarrow \infty} 
d(\vec{\Phi}_k(B_{\varepsilon}(a_i)), \vec{\Phi}_k(\partial B_{\varepsilon}(a_i)))> 0 \quad ,
\end{displaymath}
where $d$ denotes the usual distance between two sets in $\mathbb{R}^n$ \\
$d(A,B):= \sup_{p \in B} \inf_{q \in A} \vert p-q \vert$, 
then 
\begin{displaymath}
\liminf_{\varepsilon \rightarrow 0} \liminf_{k \rightarrow \infty} W(\vec{\Phi}_k (B_{\varepsilon}(a_i))) \geq 4 \pi.
\end{displaymath}
\item[ii)]
Assume that there exists an index $i$ such that 
\begin{displaymath}
\liminf_{\varepsilon \rightarrow 0} \liminf_{k \rightarrow \infty} 
Area(\vec{\Phi}_k(B_{\varepsilon}(a_i)) )>0
\end{displaymath}
then 
\begin{displaymath}
\liminf_{\varepsilon \rightarrow 0} \liminf_{k \rightarrow \infty} W(\vec{\Phi}_k (B_{\varepsilon}(a_i))) \geq 4 \pi.
\end{displaymath}
\end{itemize}
\end{Lemma}

\textit{Proof:}\\
Before passing to the proof of $i)$ let us start with a general observation (which holds independently of $i)$ or $ii)$ above): we claim that, up to passing to a subsequence in $k$, we have 
\begin{equation}\label{eq:H1to0}
\lim_{\varepsilon \rightarrow 0}\lim_{k \rightarrow \infty} \mathcal{H}^1(\vec{\Phi}_{k}(\partial B_{\varepsilon}(a_i))) \rightarrow 0.
\end{equation}
Indeed, recall that we have
\begin{displaymath}
\vec{\Phi}_k \rightarrow \vec{\Phi}_{\infty} \; \textrm{weakly in} \; W^{2,2}(\Sigma^2 \backslash \cup_{i=1}^n B_{\varepsilon}(a_i))\quad .
\end{displaymath}
Then from the classical trace theorem (see for instance \cite{RuSi}), we know that $\nabla\vec{\Phi}_k(\partial B_{\varepsilon}(a_i))$ converges weakly in $F^{\frac{1}{2}}_{2,2}(\partial B_{\varepsilon}(a_i))=H^{\frac{1}{2}}(\partial B_{\varepsilon}(a_i))$. Here,  $F^s_{p,q}$ denotes the standard Triebel-Lizorkin space. Therefore, from the standard compactness part of the Sobolev embedding theorem, we find that $\nabla\vec{\Phi}_k(\partial B_{\varepsilon}(a_i))$ converges even strongly in $L^2(\partial B_{\varepsilon}(a_i))$ and in particular in $L^1(\partial B_{\varepsilon}(a_i))$  . \\
Thus we have
\begin{displaymath}
\mathcal{H}^1(\vec{\Phi}_k (\partial B_{\varepsilon}(a_i))) = 
\int_{\partial B_{\varepsilon}(a_i)} \vert \dot{\vec{\Phi}}_k \vert \, dl \rightarrow \mathcal{H}^1(\vec{\Phi}_{\infty} (\partial B_{\varepsilon}(a_i))) = \int_{\partial B_{\varepsilon}(a_i)} \vert \dot{\vec{\Phi}}_{\infty} \vert \, dl .
\end{displaymath}
Finally, we recall that the limit immersion $\vec{\Phi}_{\infty}$ is Lipschitz, so 
\begin{displaymath}
\lim_{\varepsilon \to 0}\int_{\partial B_{\varepsilon}(a_i)} \vert \dot{\vec{\Phi}}_{\infty} \vert \, dl \rightarrow 0\quad ,
\end{displaymath}
since $\vert \partial B_{\varepsilon}(a_i) \vert \rightarrow 0$ as $\varepsilon \rightarrow 0$. \\

\textit{Proof of part i):} \\
In order to prove statement $i)$, we will use a result of Rivi\`ere (see \cite{Riv4}) giving an estimate on the Willmore energy of a compact surface with boundary; for the precise statement see Lemma \ref{LemmaRiv1} in the Appendix. Applying the  mentioned estimate to the immersions $\vec{\Phi}_k$ restricted to $B_\varepsilon(a_i)$, we find
\begin{eqnarray}
4 \pi &\leq& \liminf_{\varepsilon \rightarrow 0} \liminf_{k \rightarrow \infty} \Big(
W(\vec{\Phi}_k(B_{\varepsilon}(a_i))) +2\frac{\mathcal{H}^1(\vec{\Phi}_k(\partial B_{\varepsilon}(a_i)))}{d(\vec{\Phi}_k(\partial B_{\varepsilon}(a_i)), \vec{\Phi}_k(B_{\varepsilon}(a_i)))} \Big) \nonumber \\
&=& \liminf_{\varepsilon \rightarrow 0} \liminf_{k \rightarrow \infty} 
W(\vec{\Phi}_k(B_{\varepsilon}(a_i))) \quad,
\end{eqnarray}
where in the last equality we used \eqref{eq:H1to0} together the  assumption of $i)$. Thus, assertion $i)$ is proved.\\

\textit{Proof of part ii):}
Without loss 
of generality we might assume that the hypothesis holds for the point $a_1$. 
In the sequel it is enough to study the situation around this particular point.\\

First of all, let us recall the following assertion which can be seen as a version of Green's theorem for immersed surfaces: \\
Let $\Sigma^2$ be a smooth compact surface with boundary, 
let $\vec{\Phi} \in \mathcal{E}_{\Sigma^2}$ be an immersion of $\Sigma^2$ into $\mathbb{R}^3$
and let $\vec{X}$ be a smooth vector field in $\mathbb{R}^3$. Then it holds
\begin{equation}\label{Gauss}
\int_{\vec{\Phi}(\Sigma^2)} div_{\vec{\Phi}(\Sigma^2)} \vec{X}\, dvol_{g_{\vec{\Phi}}} =
\int_{\vec{\Phi}(\partial \Sigma^2)} \langle \vec{X}, \vec{\nu} \rangle \, dl_{\vec{\Phi}(\partial\Sigma^2)}
-2 \int_{\vec{\Phi}(\Sigma^2)} \langle \vec{X}, \vec{H} \rangle \, dvol_{g_{\vec{\Phi}}}
\end{equation}
where 
\begin{displaymath}
div_{\vec{\Phi}(\Sigma^2)} \vec{X} := \sum_{k=1}^2 \langle d\vec{X} \cdot \vec{e}_k , \vec{e}_k \rangle
\end{displaymath}
for any local orthonormal frame $(\vec{e}_1, \vec{e}_2)$ on $\vec{\Phi}(\Sigma^2)$ and 
$\vec{\nu}:= e^{-\lambda}\partial_r \vec{\Phi}$ (where as usual $\lambda$ denotes the conformal 
factor of the immersion $\vec{\Phi}$.\\
A proof of this  classical formula can be found for instance in \cite{Riv4}.\\

Now, we continue the proof of part $ii)$ of the Lemma \ref{lemmaBubbling}. \\
In a first step, we will apply formula (\ref{Gauss}) to $\vec{\Phi}_k$, restricted to the
ball $B_{\varepsilon}(a_1)$, and  to the vector field $\vec{X}(x)=x-\vec{\Phi}_k(a_1)$. Observing that $div X \equiv 2$, we get
\begin{equation}
\begin{split}
&2Area(\vec{\Phi}_k(B_{\varepsilon}(a_1))) =
\int_{\vec{\Phi}_k(B_{\varepsilon}(a_1))} div \vec{X}\, dvol_{g_k} \nonumber \\
&\quad =\int_{\vec{\Phi}_k(\partial B_{\varepsilon}(a_1))} \langle \vec{X}, \vec{\nu} \rangle \, dl
-2 \int_{\vec{\Phi}_k(B_{\varepsilon}(a_1))} \langle \vec{X}, \vec{H} \rangle \, dvol_{g_k} \nonumber \\
&\quad \leq diam(\vec{\Phi}_k(B_{\varepsilon}(a_1))) \;  \mathcal{H}^1(\partial \vec{\Phi}_k(B_{\varepsilon}(a_1)))
-2 \int_{\vec{\Phi}_k(B_{\varepsilon}(a_1))} \langle \vec{X}, \vec{H} \rangle \, dvol_{g_k} \nonumber \\
&\quad \leq diam(\vec{\Phi}_k(B_{\varepsilon}(a_1))) \; \mathcal{H}^1(\partial \vec{\Phi}_k(B_{\varepsilon}(a_1))) \nonumber \\
&\quad \quad +2 diam(\vec{\Phi}_k(B_{\varepsilon}(a_1))) \; W(\vec{\Phi}_k(B_{\varepsilon}(a_1)))^{\frac{1}{2}} \; Area(\vec{\Phi}_k(B_{\varepsilon}(a_1)))^{\frac{1}{2}}. \nonumber
\end{split}
\end{equation}
Rearranging terms, we find

\begin{equation}
Area(\vec{\Phi}_k(B_{\varepsilon}(a_1)))^{\frac{1}{2}} - \frac{diam(\vec{\Phi}_k(B_{\varepsilon}(a_1))) \mathcal{H}^1(\partial \vec{\Phi}_k(B_{\varepsilon}(a_1)))}{2 Area(\vec{\Phi}_k(B_{\varepsilon}(a_1)))^{\frac{1}{2}}} \leq
diam(\vec{\Phi}_k(B_{\varepsilon}(a_1)))  \,W(\vec{\Phi}_k(B_{\varepsilon}(a_1)))^{\frac{1}{2}}.
\end{equation}
Next,  we recall that the diameters of the whole surfaces
$\vec{\Phi}_k(\Sigma^2)$ as well as the diameter of the limit $\vec{\Phi}_\infty(\Sigma^2)$ are uniformly bounded thanks to Lemma \ref{LemmaSimon}. Thus, we estimate 
\begin{equation}
\begin{split}
&\sqrt{8\pi} \; diam(\vec{\Phi}_k(B_{\varepsilon}(a_1)))  \quad \geq W(\vec{\Phi}_k(B_{\varepsilon}(a_1))^{\frac{1}{2}} \;  diam(\vec{\Phi}_k(B_{\varepsilon}(a_1))) \nonumber \\
& \quad \geq 
Area(\vec{\Phi}_k(B_{\varepsilon}(a_1)))^{\frac{1}{2}} - \frac{diam(\vec{\Phi}_k(B_{\varepsilon}(a_1))) \; \mathcal{H}^1(\partial \vec{\Phi}_k(B_{\varepsilon}(a_1)))}{2 Area(\vec{\Phi}_k(B_{\varepsilon}(a_1)))^{\frac{1}{2}}} \nonumber \\
& \quad \geq \eta^{\frac{1}{2}}
- \frac{diam(\vec{\Phi}_k(B_{\varepsilon}(a_1))) \; \mathcal{H}^1(\partial \vec{\Phi}_k(B_{\varepsilon}(a_1)))}{2 Area(\vec{\Phi}_k(B_{\varepsilon}(a_1)))^{\frac{1}{2}}} \nonumber \\
& \quad \geq \eta^{\frac{1}{2}} 
- \frac{C \, \mathcal{H}^1(\partial \vec{\Phi}_k(B_{\varepsilon}(a_1)))}{2 Area(\vec{\Phi}_k(B_{\varepsilon}(a_1)))^{\frac{1}{2}}}, \quad \text{for some $C>0$} \nonumber \\
& \quad > \eta_1 >0 ,\quad  \text{since $\mathcal{H}^1(\partial \vec{\Phi}_k(B_{\varepsilon}(a_1))) \rightarrow 0.$} \nonumber
\end{split}
\end{equation}
Therefore  we found that $diam(\vec{\Phi}_k(B_{\varepsilon}(a_1)))\geq \eta_2 >0$ for some 
$\eta_2$. But this  together with the fact that $\mathcal{H}^1(\partial \vec{\Phi}_k(B_{\varepsilon}(a_1))) \rightarrow 0$
implies that  $$d(\vec{\Phi}_k(\partial B_{\varepsilon}(a_1)), \vec{\Phi}(B_{\varepsilon}(a_1)))\geq \eta^*>0$$
for some positive $\eta^*$. We can now conclude the proof of part $ii)$ by applying the (already proved) part $i)$. 
\hfill$\Box$\\

Now, let us explain how Lemma \ref{lemmaBubbling} implies that the limit satisfies the isoperimetric constraint. We claim that, for every  concentration point $a_i$,
\begin{equation}\label{eq:claim1}
\liminf_{\varepsilon \rightarrow 0} \liminf_{k \rightarrow \infty} 
d(\vec{\Phi}_k(B_{\varepsilon}(a_i)), \vec{\Phi}_k(\partial B_{\varepsilon}(a_i)) = 0
\end{equation}
and
\begin{equation}\label{eq:claim2}
\liminf_{\varepsilon \rightarrow 0} \liminf_{k \rightarrow \infty} 
Area(\vec{\Phi}_k(B_{\varepsilon}(a_i)) =0.
\end{equation}
Indeed, if by contradiction there exists a concentration point, say $a_1$, where one of the two statements above fail then   
\begin{eqnarray}
\liminf_{k \rightarrow \infty}W(\vec{\Phi}_k)
&=& \lim_{\varepsilon \rightarrow 0 }\liminf_{k \rightarrow \infty} \left[ \big(W(\vec{\Phi}_k(\Sigma^2\backslash \cup_{i=1}^n B_{\varepsilon}(a_i))) 
+ W(\vec{\Phi}_k(\cup_{i=1}^n B_{\varepsilon}(a_i))) \right]\nonumber \\
&\geq& W(\vec{\Phi}_\infty) + 4\pi \geq 8\pi \quad, 
\end{eqnarray}
where the first inequality comes from the lower semicontinuity of $W$ under weak-$W^{2,2}$ convergence together with the fact that $\vec{\Phi}_\infty$ is an element of ${\mathcal E}_{\Sigma^2}$ and Lemma \ref{WillmoreInf} applied to $a_1$; the last inequality follows from $W(\vec{\Phi}_\infty)\geq 4\pi$. 
\\

But this last estimate is in contradiction with our hypothesis that the Willmore energy stays strictly below $8 \pi$. Thus, our claim \eqref{eq:claim1}-\eqref{eq:claim2} holds. 
\\

Now, combining Remark \ref{rem:IsopLimit} and the second claim  \eqref{eq:claim2}  we deduce that the ares converge, i.e.
\begin{displaymath}
Area(\vec{\Phi}_k)\rightarrow Area(\vec{\Phi}_{\infty}).
\end{displaymath}
Moreover, combining \eqref{eq:H1to0}  and the first claim \eqref{eq:claim1}, we infer that  also the diameters have to vanish:
\begin{displaymath}
diam(\vec{\Phi}_k(B_{\varepsilon}(a_i))) \rightarrow 0 \; \textrm{for all} \, i.
\end{displaymath}
But this last fact together with Remark \ref{rem:IsopLimit} immediately leads to the conclusion that the enclosed volumes converge as well
\begin{displaymath}
Vol\Big( \; \overline{\vec{\Phi}}_k \;\Big) \rightarrow Vol\Big( \; \overline{\vec{\Phi}}_{\infty}\;\Big).
\end{displaymath}
We conclude that the isoperimetric constraint passes to the limit.

\subsubsection{Regularity of the limit, proof of Corollary \ref{Corollary}}\label{Subsub:reg}

First of all, thanks to Corollary \ref{cor:embedding} observe that the uniform energy estimate
\begin{displaymath}
W(\vec{\Phi}_k) < 8 \pi - \delta
\end{displaymath}
together with the lower semicontinuity of the Willmore functional under weak $W^{2,2}$ convergence immediately implies that the limit $\vec{\Phi}_{\infty}$ is an embedding, i.e. $\vec{\Phi}_\infty$ has no self intersection. Moreover, as we have seen in the preceding  subsection, the limit satisfies the isoperimetric constraint.\\

Now, we  study the regularity question.  In a first step we will recall that the Willmore functional $W$ is Fr\'echet differentiable and we will determine $dW$. Then we will establish the equation which is satisfied by the limit $\vec{\Phi}_{\infty}$ and finally show the regularity.

\begin{Lemma}
Let $\vec{\Phi}$ belong to $\mathcal{E}_{\Sigma^2}$. Then the Willmore functional $W$ is Fr\'echet differentiable with respect to variations $\vec{w} \in W^{1,\infty} \cap W^{2,2}$ with compact support. \\
Moreover, for the differential we have the following formula
\begin{displaymath}
d_{\vec{\Phi}}W[\vec{w}]= \int \nabla \vec{w} \cdot \Big ( 2\nabla \vec{H} -3 H \nabla \vec{n} - \nabla^{\perp} \vec{n} \times \vec{H}  \Big ).
\end{displaymath}
\end{Lemma}

This lemma is a straightforward adaptation of the corresponding result in \cite{Riv5} (see also  the Appendix of \cite{MR2}).
Note that the differential exactly corresponds to the reformulation of the Willmore equation in the form of a conservation law as in (\ref{EqWillmore}).\\

Next, we give the equation which is satisfied by our limiting object $\vec{\Phi}_{\infty}$.
\begin{Lemma}
Let $\vec{\Phi}_{\infty}$ be as above. Then it satisfies the following equation
\begin{displaymath}
 div ( 2\nabla \vec{H} -3 H \nabla \vec{n} - \nabla^{\perp} \vec{n} \times \vec{H}- \mu \nabla \vec{\Phi} )=0.
\end{displaymath}
\end{Lemma}

\textit{Proof:}  The assertion is an immediate consequence of the lower semicontinuity the  Willmore functional $W$ and its differentiability, the classical fact that the area functional $A$ is Fr\'echet differentiable as well with differential
\begin{displaymath}
d_{\vec{\Phi}}A[\vec{w}]= - \int 2 \vec{H} \cdot \vec{w} 
\end{displaymath}
and the principle of Lagrange multipliers. To get the final formula recall also that $ \vec{H}=\frac{1}{2} \Delta_g \vec{\Phi} = \frac{1}{2}e^{-2\lambda} \Delta \vec{\Phi}$. \hfill$\Box$\\

Finally we have the following regularity assertion, whose proof is an easy adaptation of the regularity of Willmore immersions established in \cite{Riv2} (for the regularity of area-constrained Willmore immersions see also \cite{MR2}, and for a comprehensive explanation see \cite{RivNotes}).
\begin{Lemma}
$\vec{\Phi}_{\infty}$ is smooth. 
\end{Lemma}

The proof of Corollary \ref{Corollary} is now complete. \hfill$\Box$

\subsection{The case of diverging conformal factors}\label{Subsec:Diverging}

So far, we have studied the case when $C_k$ remain bounded. Now, we will analyze the case where this is no longer true. As we have seen 
above, the only remaining possibility is that - up to extraction of a subsequence - the 
conformal factors tend to $-\infty$.\\

\subsubsection{Existence of at most one bubble}\label{Subsub:OneBubble}

As above, the points where Willmore energy concentrates are denoted by $a_i$, and again we may assume that we have
\begin{displaymath}
Area(\vec{\Phi}_k) = 1 \; \forall \, k.
\end{displaymath}
Thus, the fact that now the conformal factors diverge, i.e.
\begin{displaymath}
\lambda_k \rightarrow - \infty \; \textrm{on} \, \Sigma^2 \backslash \cup_{i=1}^n B_{\varepsilon}(a_i)
\end{displaymath}
implies that there is at least one point $a_i$ where area concentrates.\\ 
More precisely, we have that there exists a point $a^*$ such that
\begin{equation}\label{AreaEstimate}
\liminf_{\varepsilon \rightarrow 0} \liminf_{k \rightarrow \infty} 
Area(\vec{\Phi}_k(B_{\varepsilon}(a_i))) > 0.
\end{equation}

Our next goal is to show that there is exactly one such point; this is the content of the following lemma, whose proof follows by the $8\pi-\delta$ bound on the Willmore energy of $\vec{\Phi}_k$ and  the second part of Lemma \ref{lemmaBubbling}.

\begin{Lemma}\label{Lemma1Point}
Assume that the conformal factors tend to $-\infty$, i.e.
\begin{displaymath}
\lambda_k \rightarrow -\infty \; \textrm{on} \; \Sigma^2 \backslash \cup_{i=1}^n B_{\varepsilon}(a_i)\; \textrm{as} \, k \rightarrow \infty.
\end{displaymath}
Then there exists exactly one point $a^* \in \Sigma^2$  of concentration of the area, i.e. where  where \eqref{AreaEstimate} holds. Of course $a^*$ is also a point of concentration of the energy, i.e. $a^*\in \{a_i\}$.  
\end{Lemma}

Directly from the lemma it follows that
\begin{eqnarray}
\liminf_{\varepsilon \rightarrow 0} \liminf_{k \rightarrow \infty}
Area(\vec{\Phi}_k(B_{\varepsilon}(a^*))) =1  \quad \text{and } \label {eq:Area=1} \\ 
\limsup_{\varepsilon \rightarrow 0} \limsup_{k \rightarrow \infty}
Area(\vec{\Phi}_k(\Sigma^2 \backslash B_{\varepsilon}(a^*))) =0  \label{eq:Areato0} \quad .
\end{eqnarray}

Analogously, again using the $8\pi-\delta$ bound on the Willmore energy of $\vec{\Phi}_k$ and  the first part of Lemma \ref{lemmaBubbling}, we also have
\begin{eqnarray}
\liminf_{\varepsilon \rightarrow 0} \liminf_{k \rightarrow \infty} 
d(\vec{\Phi}_k(B_{\varepsilon}(a^*)), \vec{\Phi}_k(\partial B_{\varepsilon}(a_i)) >0 \quad \text{and } \label{eq:d>0} \\
\liminf_{\varepsilon \rightarrow 0} \liminf_{k \rightarrow \infty} d(\vec{\Phi}_k(\Sigma^2 \setminus B_{\varepsilon}(a^*)), \vec{\Phi}_k(\partial B_{\varepsilon}(a_i)) =  0 \label{eq:dto0}.
\end{eqnarray}
On the other hand, recalling \eqref{eq:H1to0}, also the lengths of the boundaries converge to zero
\begin{displaymath}
\limsup_{\varepsilon \rightarrow 0}\limsup_{k \rightarrow \infty} \mathcal{H}^1(\vec{\Phi}_{k}(\partial B_{\varepsilon}(a^*))) \rightarrow 0,
\end{displaymath}
and combining this last observation with \eqref{eq:dto0} we obtain that also the diameter of the complementary of the bubble has to vanish in the limit:
\begin{equation}\label{eq:diamBcto0}
\limsup_{\varepsilon \to 0} \limsup_{k\to \infty} diam(\vec{\Phi}_k(\Sigma^2\setminus B_\varepsilon(a^*))) =0 \quad .
\end{equation}

Observe that \eqref{eq:Areato0}, respectively \eqref{eq:diamBcto0}, implies that the portion of the surface outside the bubble in $a^*$ is not contributing in the limit to the area, respectively to the enclosed volume; therefore it does not contribute to the isoperimetric ratio in the limit. Nevertheless notice that the bubble forming in $a^*$ is topologically (a bigger and bigger portion of) a  sphere, therefore the topological information is carried by the shrinking part made by the complementary of the bubble. Summarizing, we are facing a \emph{dichotomy} between a (bigger) portion of the surface (namely, the bubble) carrying the isoperimetric ratio, and a (smaller) portion (the complementary of the bubble) carrying the topological information.  
\\In the next subsection we isolate the two parts by performing  a ``cut and fill'' procedure.  

\subsubsection{Cut and fill}
 
As described above the geometric situation can be described as follows (for $k$ sufficiently large): there is a dichotomy between the topological information and the additional constraint of prescribed isoperimetric ratio. More precisely, $\vec{\Phi}_k(B_{\varepsilon}(a^*))$ forms a spherical bubble
carrying the isoperimetric information and $\vec{\Phi}_k(\Sigma^2 \backslash B_{\varepsilon}(a^*))$ keeps  the topological information, i.e. 
it has the topological type of $\Sigma^2$, i.e. a genus $\g\geq 1$ surface. \\ 

The strategy now is to find estimates for the Willmore energy for the two parts $\vec{\Phi}_k(\Sigma^2 \backslash B_{\varepsilon}(a^*))$ 
and $\vec{\Phi}_k(B_{\varepsilon}(a^*))$ and bring these estimates to a contradiction to our additional hypothesis relating the Willmore energy of Schygulla-spheres to 
the Willmore energy of our embeddings $\vec{\Phi}_k$. \\
In doing so, we will exploit the existence of a genus $\g\geq 1$ minimizers of the Willmore energy among all genus $\g$ embedded surfaces (free minimization), the existence of Schygulla-spheres $\mathbb{S}_{S,r}$, i.e. smoothly embedded surfaces of type $\mathbb{S}^2$  minimizing the Willmore energy for the given isoperimetric ratio $r$, and the fact
that the function which assigns to each given isoperimetric ratio the Willmore energy of the Schygulla-sphere for this given isoperimetric ratio is strictly monotone and continuous. These latter facts are proved
in \cite{Schy}.\\

In order to perform the above strategy, we will apply a cut-and-fill-procedure in order to close each of the parts $\vec{\Phi}_k(\Sigma^2 \backslash B_{\varepsilon}(a^*))$ and $\vec{\Phi}_k(B_{\varepsilon}(a^*))$, complete them separately to new closed surfaces and finally estimate appropriately their Willmore energy. 
For this procedure, we will use the following adaptation of a lemma which can be found (with proof)
in \cite{MR}. 

\begin{Lemma}\label{CutAndFill}
Let $\vec{\Phi}_k$ be a sequence of conformal weak, 
immersions $\{\vec{\Phi}_k\}_{k \in \N} \subset \mathcal{E}_{\Sigma^2}$ into $\mathbb{R}^3$.
Assume that
\begin{displaymath}
\limsup_{k \rightarrow \infty} \int_{\Sigma^2} \lbrack 1 + \vert D\vec{n}_{\vec{\Phi}_k} \vert^2 \rbrack dvol_{g_k} < \infty.
\end{displaymath}
Let $a \in \Sigma^2$ and $s_k$, $t_k \rightarrow 0$ such that
\begin{displaymath}
\frac{t_k}{s_k} \rightarrow 0
\end{displaymath}
and
\begin{displaymath}
\lim_{k \rightarrow \infty} \int_{B_{s_k}(a) \backslash B_{t_k}(a)} \lbrack 1 + \vert \mathbb{I}_{\vec{\Phi}_k}\vert^2 \rbrack dvol_{g_k} =0.
\end{displaymath}
Then there exist conformal immersions $\vec{\xi}_k$ from $\Sigma^2$ into $\mathbb{R}^3$ and a sequence of quasi conformal bilipschitz
homeomorphisms 	$\psi_k$ of $\Sigma^2$, converging in $C^0$-norm over $\Sigma^2$ to the identity map, such that
\begin{displaymath}
\vec{\xi}_k \circ \psi_k = \vec{\Phi}_k \; \textrm{in} \; \Sigma^2 \backslash B_{s_k}(a)
\end{displaymath}
and
\begin{displaymath}
\lim_{k \rightarrow \infty} diam(\vec{\xi}_k \circ \psi_k (B_{s_k}(a))=0, \; \lim_{k \rightarrow \infty} Area(\vec{\xi}_k\circ \psi_k (B_{s_k}(a))=0.
\end{displaymath}
Moreover
\begin{displaymath}
\lim_{k \rightarrow \infty} \int_{B_{s_k}(a)}  \vert \mathbb{I}^0_{\vec{\xi}_k \circ \psi_k} \vert^2  dvol_{g_{\vec{\xi}_k \circ \psi_k}} =0
\end{displaymath}
where $\mathbb{I}^0_{\vec{\xi}_k \circ \psi_k}$ is the trace free second fundamental form.
\end{Lemma}

\begin{Remark}\label{rem:CutFill}
\rm{
\begin{itemize}
\item
Note that due to the assumption of bounded Willmore energy
\begin{displaymath}
\inf_k W(\vec{\Phi}_k)< 8\pi
\end{displaymath}
the hypothesis of the above lemma are satisfied, in particular  it is possible to find radii $s_k$ and $t_k$ 
such that
\begin{displaymath}
\lim_{k \rightarrow \infty} \int_{B_{s_k}(a^*) \backslash B_{t_k}(a^*)} \lbrack 1 + \vert \mathbb{I}_{\vec{\Phi}_k}\vert^2 \rbrack dvol_{g_k} =0.
\end{displaymath}
\item
At first glimpse, one might have the impression that by the cut-and-fill-procedure guaranteed by the previous lemma we can close just one of the two parts of our surface, either $\vec{\Phi}_k(\Sigma^2 \backslash B_{\varepsilon}(a^*))$ or $\vec{\Phi}_k(B_{\varepsilon}(a^*))$. But the proof of the 
above cited result actually reveals that the procedure allows to close both parts. More precisely, depending on the side  from which we approach the curve along which we cut (from ``inside'' or from ``outside'') we have two possibilities of closing the surface: either we glue almost an entire shrinking sphere $\Sigma^2_{1,s_k}$ with Willmore energy $W(\Sigma^2_{1,s_k}) \rightarrow 4 \pi$ as $s_k \rightarrow 0$ or we glue a shrinking almost flat disk $\Sigma^2_{2,s_k}$ with Willmore energy $W(\Sigma^2_{2,s_k}) \rightarrow 0$ as $s_k \rightarrow 0$.
\item
The only remaining risk would be that we create branch points. \\
From the proof of the original version of the above lemma in \cite{MR} we see that this can be excluded once we can show that in 
\begin{displaymath}
\Delta \hat{\lambda}_{\infty} = c_0 \delta_0
\end{displaymath}
it holds $c_0=0$. \\
In the above equation, $\hat{\lambda}_{\infty}$ is the limit of the conformal factors 
$\hat{\lambda}_k= \log \vert \partial_{x_1} \hat{\vec{\Phi}}_k\vert= \log \vert \partial_{x_2} \hat{\vec{\Phi}}_k\vert$ where
\begin{displaymath}
\hat{\vec{\Phi}}_k:= e^{-C_k}(\vec{\Phi}_k - \vec{\Phi}_k(x_0))
\end{displaymath}
for a suitable $x_0$, and $\delta_0$ denotes the Dirac delta distribution centered at the origin. We claim that, in the present framework, actually we have $c_0=0$. \\
Indeed, by assumption, recall that there exists $\delta >0$ such that
\begin{equation}\label{eq:W<8pi}
W(\vec{\Phi}_k) \leq 8 \pi -\delta \; \textrm{for all $k$}.
\end{equation}
If we had $c_0 \neq 0$, $\hat{\vec{\Phi}}_{\infty}$ had a branch point of order $\frac{c_0}{2 \pi}$. In other words, we had that 
$\hat{\vec{\Phi}}_{\infty}$ covers $\frac{c_0}{2 \pi}+1$ times the plane $P^2_0$. Due to the fact that for any given $\varepsilon > 0$, on $\Sigma^2 \backslash B_{\varepsilon}(a^*)$ we have weak $W^{2,2}$-convergence of $\hat{\vec{\Phi}}_k$ to $\hat{\vec{\Phi}}_\infty$, then for any choice of $0< 2 \alpha < \beta < \varepsilon$
\begin{displaymath}
\hat{\vec{\Phi}}_k \rightarrow \hat{\vec{\Phi}}_{\infty} \; \textrm{in} \; C^{0, \gamma}(B_{\beta}(a^*) \backslash B_{2 \alpha}(a^*))
\end{displaymath}
by the classical Sobolev embedding theorem. \\
In this situation, we select a point $p \in B_{\beta}(a^*) \backslash B_{2 \alpha}(a^*)$ and a radius $\eta$ small enough such that
\begin{displaymath}
B_{\eta}(p) \subset B_{\beta}(a^*) \backslash B_{2 \alpha}(a^*)
\end{displaymath}
and we apply Lemma \ref{LemmaSimon2} in the Appendix
to $\vec{\Phi}_k(\Sigma^2)$ and a ball $B_{\rho} \subset \mathbb{R}^3$ such that $\vec{\Phi}_k(B_{\eta}(p) ) \subset B_{\rho}$ \ in order to conclude that - for $\eta$ and $\rho$ small enough - 
\begin{displaymath}
\limsup_k W(\vec{\Phi}_k) \geq 8 \pi - \frac{\delta}{2}.
\end{displaymath}
Since this contradicts our hypothesis \eqref{eq:W<8pi}, we must have $c_0=0$ and therefore the immersions we created via Lemma \ref{CutAndFill} are unbranched.
\end{itemize}}
\hfill$\Box$
\end{Remark}

In a first step, applying Lemma \ref{CutAndFill}, we get the following equalities for $k$ large enough
\begin{eqnarray}
W(\vec{\Phi}_k) 
&=& W(\vec{\Phi}_k(\Sigma^2 \backslash B_{s_k}(a^*))) + W(\vec{\Phi}_k (B_{s_k}(a^*)))   \nonumber \\
&=& W(\vec{\xi}_k \circ \psi_k (\Sigma^2 \backslash B_{s_k}(a^*))) + W(\vec{\Phi}_k (B_{s_k}(a^*))) .\nonumber 
\end{eqnarray}
Now, in the present dichotomy situation we have the following two possibilities of filling: either we glue an almost entire sphere $\Sigma^2_{1,s_k}$ to $\vec{\xi}_k \circ \psi_k (\Sigma^2 \backslash B_{s_k}(a^*))$ and an almost flat disc $\Sigma^2_{2,s_k}$ to $\vec{\Phi}_k (B_{s_k}(a^*))$ or viceversa. \\
In the first case we can estimate 
\begin{eqnarray}
W(\vec{\Phi}_k) &=&
W(\vec{\xi}_k \circ \psi_k (\Sigma^2 \backslash B_{s_k}(a^*))) + W(\vec{\Phi}_k (B_{s_k}(a^*)))  \nonumber \\
&=& W(\vec{\xi}_k \circ \psi_k (\Sigma^2 \backslash B_{s_k}(a^*)) \cup \Sigma^2_{1,s_k}) -W(\Sigma^2_{1,s_k})  \nonumber \\
&& \quad + W(\vec{\Phi}_k (B_{s_k}(a^*)) \cup \Sigma^2_{2,s_k}) - W(\Sigma^2_{2,s_k}) \nonumber \\
&\geq& \beta^3_{\g} -W(\Sigma^2_{1,s_k})  + W(\vec{\Phi}_k (B_{s_k}(a^*)) \cup \Sigma^2_{2,s_k}) - W(\Sigma^2_{2,s_k}) \nonumber \\
&& \quad \textrm{since, by construction, $\vec{\xi}_k \circ \psi_k (\Sigma^2 \backslash B_{s_k}(a^*)) \cup \Sigma^2_{1,s_k}$ is a genus $\g$ surface } \nonumber \\
&\geq&  \beta^3_{\g} - 4 \pi - \varepsilon_1(s_k)  + W(\vec{\Phi}_k (B_{s_k}(a^*)) \cup \Sigma^2_{2,s_k}) - W(\Sigma^2_{2,s_k}) \nonumber \\
&& \textrm{by the energy estimate for $\Sigma^2_{1,s_k}$ of Remark \ref{rem:CutFill}} \nonumber \\
&\geq&  \beta^3_{\g} - 4 \pi - \varepsilon_1(s_k) + W(\mathbb{S}_{S, iso(\vec{\Phi}_k (B_{s_k}(a^*)) \cup \Sigma^2_{2,s_k})}) - W(\Sigma^2_{2,s_k}) \nonumber \\
&& \textrm{since $\vec{\Phi}_k (B_{s_k}(a^*)) \cup \Sigma^2_{2,s_k}$ has the topology of a sphere} \nonumber \\
&\geq& \beta^3_{\g} - 4 \pi - \varepsilon_1(s_k)+ W(\mathbb{S}_{S, iso(\vec{\Phi}_k (B_{s_k}(a^*)) \cup \Sigma^2_{2,s_k})}) - \varepsilon_2(s_k) \nonumber \\
&& \textrm{by the energy estimate for $\Sigma^2_{2,s_k}$ of Remark \ref{rem:CutFill}} \quad . \nonumber 
\end{eqnarray}
Observe that as $s_k$ tend to $0$, the two terms $\varepsilon_1(s_k)$ and $\varepsilon_2(s_k)$ vanish.\\
Moreover, exploiting the fact that 
\begin{displaymath}
t \mapsto \beta(t):=W(\mathbb{S}_{S,t})
\end{displaymath} 
is continuous in $t$ (see \cite{Schy}) together with the observation that  (here we use that $\vec{\Phi}_k (B_{s_k}(a^*))$ carries the isoperimetric information of $\vec{\Phi}_k$, and the isoperimetric contribution of $\Sigma^2_{2,s_k}$ is negligible in the limit)
$$\lim_{k \to \infty} iso(\vec{\Phi}_k (B_{s_k}(a^*)) \cup \Sigma_{2,s_k}) = R \quad ,$$ 
we conclude that
\begin{displaymath}
W(\mathbb{S}_{S, iso(\vec{\Phi}_k (B_{s_k}(a^*)) \cup \Sigma^2_{2,s_k})}) \rightarrow W(\mathbb{S}_{S,R}).
\end{displaymath}
In the second case (i.e. the shrinking almost entire sphere $\Sigma^2_{1,s_k}$ is glued to the bubble $\vec{\Phi}_k(B_{s_k}(a^*))$, and the almost flat disc $\Sigma^2_{2,s_k}$ is glued to $\vec{\xi}_k \circ \psi_k (\Sigma^2 \backslash B_{s_k}(a^*))$ ) the analogous estimate bring to the same conclusions.\\

So, in the limit we find
\begin{displaymath}
\liminf_ {k \rightarrow \infty}W(\vec{\Phi}_k) \geq W(T_{Clifford}) + W(\mathbb{S}_{S,R}) - 4 \pi,
\end{displaymath}
which contradicts that $R\in I_{\g}$, and $I_{\g}$ is defined as in \eqref{eq:defI}. Therefore the case $C_k \to -\infty$ cannot occur. Since the only remaining case is when $\{C_k\}_{k \in \N}$ is bounded, and in Subsection \ref{Subsec:CkBounded} we already proved the existence of a minimizer under this assumption, the proof of Theorem \ref{MainTheorem} is now complete.  \hfill$\Box$

\section{Proof of Theorem \ref{Thm2}}

First of all, we will show that the set $I_{\g}$ is not empty. To this aim recall that for every genus $\g\geq 1$ the infimum $\beta^3_{\g}$ of the Willmore energy among genus $\g$ immersed surfaces in $\mathbb{R}^3$ is attained by a smooth embedding $\vec{\Phi}_{\g}$ (for genus one see \cite{Si}, for higher genus see \cite{BK}; for a different proof of the general case see \cite{Riv5}). We claim that $iso(\vec{\Phi}_{\g}) \in I_{\g}$, where $I_{\g}$ is defined in \eqref{eq:defI}.

By a direct comparison with the stereographic projections of the Lawson minimal surfaces of genus $\g$ in $\mathbb{S}^3$ (see e.g. the introduction of \cite{Si}),  for every $g\geq 1$ one has 
\begin{equation}\label{eq:beta<8pi}
\beta^3_{\g}<8\pi \quad;
\end{equation}
moreover, by the work of Bauer and Kuwert \cite{BK} we know that 
\begin{equation}\label{eq:beta<omega}
\beta^3_{\g} < \omega^3_{\g} \quad ;
\end{equation}
Finally, by the isoperimetric inequality in $\R^3$, the isoperimetric ratio of any embedded surface in $\R^3$ different from a round sphere is strictly larger then $iso(\mathbb{S}^2)$, then
\begin{equation}\label{eq:iso>}
iso(\vec{\Phi}_{\g}) > iso(\mathbb{S}^2) \quad.
\end{equation}
By the fact that for any Schygulla sphere $\mathbb{S}_{S,t}$  different from the round sphere we have $W(\mathbb{S}_{S,t}) > 4 \pi$ (see \cite{Schy}), we then deduce
\begin{equation*}
W(\vec{\Phi}_{\g})= \beta^3_{\g} < \beta^3_{\g} - 4\pi + W(\mathbb{S}_{S,r=iso(\vec{\Phi}_{\g})} )\quad,
\end{equation*}
which, combined with \eqref{eq:beta<8pi} and \eqref{eq:beta<omega} yields
$$W(\vec{\Phi}_{\g})< \min\{ 8\pi, \omega^3_{\g}, \beta^3_{\g} - 4\pi + W(\mathbb{S}_{S,r=iso(\vec{\Phi}_{\g})} ) \}\quad. $$
It follows that $iso(\vec{\Phi}_{\g}) \in I_{\g}$ as desired.\\

Before we continue the proof of Theorem \ref{Thm2}, let us state and prove an easy  lemma which we will use later on.
\begin{Lemma}\label{LemmaThm2}
Look at the following map
\begin{eqnarray}
\Psi : \left \{ \Sigma^2 \, : \, \Sigma^2 \, \textrm{is a smoothly embedded surface in} \, \mathbb{R}^3 \right \} &\rightarrow& \mathbb{R}\nonumber \\
\Sigma^2 &\mapsto& iso(\Sigma^2). \nonumber
\end{eqnarray} 
Then $\Sigma^2$ is a critical point of $\Psi$ if and only if $\Sigma^2$ is a round sphere.
\end{Lemma}

\textit{Proof.}
Due to the fact that the isoperimetric ratio is invariant under rescaling, we may - without loss of generality - assume that the volume enclosed 
by $\Sigma^2$ is equal to 1.\\
Thus, $\Sigma^2$ is a critical point of $\Psi$ if and only if it is a critical point of the area functional under the constraint of fixed volume. \\
But it is well known that this leads to the conclusion that $\Sigma^2$ has to be surface of constant mean curvature. \\
And finally by a famous theorem due to H. Hopf and Alexandrov (see for instance \cite{K}, it follows that $\Sigma^2$ is a round sphere. Thus the proof of
Lemma \ref{LemmaThm2} is complete.
\hfill$\Box$
\\

Now we come back to the proof of Theorem \ref{Thm2}: given $r\in I_{\g}$, we have to show that there exists $\delta>0$ such that $(r-\delta,r+\delta)\subset I_{\g}$.\\
By Theorem \ref{MainTheorem} we know that there exists a smooth embedding $\vec{\Phi}_r: \Sigma^2 \rightarrow \mathbb{R}^3$ minimizing the Willmore 
energy and respecting the given isoperimetric constraint $iso(\vec{\Phi}_r)=r$. \\
By Lemma \ref{LemmaThm2} we deduce that $\vec{\Phi}_r$ is not a critical point of $\Psi$, in other words, the differential of $\Psi$ at the point
$\vec{\Phi}_r$ is surjective. Therefore, there exists a smooth 1-parameter family $\left \{ \vec{\Phi}_s \right \}_{s \in (r-\varepsilon_0, r+ \varepsilon_0)}$ of embeddings - obtained by perturbing $\vec{\Phi}_r$ - such that for every $\varepsilon\in(0, \varepsilon_0)$ there exists $\delta > 0$ such that
\begin{equation}\label{Equation}
\Psi \Big(\left \{ \vec{\Phi}_s \right \}_{s \in (r-\varepsilon, r+ \varepsilon)} \Big) \supset (r-\delta, r + \delta).
\end{equation}
Now, recall that for the embedding $\vec{\Phi}_r$ we started with, we have the following estimate
\begin{displaymath}
W(\vec{\Phi}_r)<
\min \left \{ 8 \pi, \omega^3_{\g}, \beta^3_{\g} + W(\mathbb{S}_{S,r}) - 4 \pi \right \}.
\end{displaymath}
Thus, there exists $\eta >0$ such that
\begin{displaymath}
W(\vec{\Phi}_r) \leq
\min \left \{ 8 \pi, \omega^3_{\g}, \beta^3_{\g}  + W(\mathbb{S}_{S,r}) - 4 \pi \right \} -\eta.
\end{displaymath}
Hence, from the smoothness in $s$ of the family $\left \{ \vec{\Phi}_s \right \}_{s \in (r-\varepsilon, r+ \varepsilon)}$, for $s$ close enough to $r$ we have 
\begin{equation}\label{eq:PfThmI}
W(\vec{\Phi}_s) \leq
\min \left \{ 8 \pi, \omega^3_{\g}, \beta^3_{\g}  + W(\mathbb{S}_{S,r}) - 4 \pi \right \} -\frac{\eta}{2}.
\end{equation}
In a last step, we  exploit the continuity of the function which to a given isoperimetric ratio assigns the Willmore energy of the corresponding Schygulla sphere, more precisely the  function $t \mapsto \beta(t):=W(\mathbb{S}_{S,t})$ is continuous in $t$ (see \cite{Schy} for the proof). We deduce  that for $s$ close enough to $r$ we have
\begin{equation}\label{eq:PfThmI2}
\vert W(\mathbb{S}_{S,\Psi (\vec{\Phi}_s)}) - W(\mathbb{S}_{S,r}) \vert \leq \frac{\eta}{4}.
\end{equation}
Combining \eqref{eq:PfThmI} and \eqref{eq:PfThmI2} we get
\begin{eqnarray}
W(\vec{\Phi}_s) &\leq& \beta^3_{\g} + W(\mathbb{S}_{S,\Psi (\vec{\Phi}_s)}) - 4 \pi  -\frac{\eta}{4}  \nonumber \\
&<& \beta^3_{\g} + W(\mathbb{S}_{S,\Psi (\vec{\Phi}_s)}) - 4 \pi. \nonumber
\end{eqnarray}
This last inequality together with \eqref{eq:PfThmI} and  \eqref{Equation} conclude the proof of Theorem \ref{Thm2}.
\hfill$\Box$

\section*{Appendix}

In this section we recall some useful  results used in the main text.\\

We start with a result which relates the non-compactness of the conformal classes to an estimate for the Willmore energy from below. \\
This result can be found in  \cite{KS5} (see also  \cite{KL} and \cite{Riv4} for the higher codimentional case). 

\begin{theorem} \label{ConfClass}
Let $(\Sigma^2,c_k)$ be a sequence of closed Riemannian surfaces of genus $\g$ and conformal classes $c_k$. Assume that $\lbrack c_k \rbrack $ is diverging to the boundary of the Moduli Space of the conformal classes on  $\Sigma^2$. 
Let $\{\vec{\Phi}_k\}_{k \in \N}\subset \mathcal{E}_{\Sigma^2}$ be a sequence of weak immersions, conformal  with respect to $c_k$; then
\begin{displaymath}
\liminf_{k \rightarrow \infty} \int_{\Sigma^2}\vert \vec{H}_{\vec{\Phi}_k}\vert^2 dvol_{\vec{\Phi}_k^*g_{\mathbb{R}^3}} \geq \min \{ 8 \pi, \omega^3_{\g} \}
\end{displaymath}
where $\omega^3_{\g}$ was defined in \eqref{eq:defomega}.
\end{theorem}

We continue with the famous Wente estimate. 

\begin{theorem} \label{WenteEtAl}
Let $a$ and $b$ be two function in $W^{1,2}(D^2, \mathbb{R})$. Moreover, let $\phi$ be the unique solution of 
\begin{displaymath}
\left\{
\begin{array}{rll}
-\Delta \varphi &= \nabla a \cdot \nabla^{\perp}b = \partial_x a \partial_y b - \partial_y a \partial_x b &\textrm{in} \; B^2_1(0)\\
\varphi&=0 & \textrm{on} \; \partial B^2_1(0),
\end{array} \right.
\end{displaymath}
Then the following estimates hold 
\begin{displaymath}
\vert \vert \phi \vert \vert_{\infty} + \vert \vert \nabla \phi \vert \vert _{2,1} + \vert \vert \nabla^2 \phi \vert \vert_1 \leq C\vert \vert \nabla a \vert \vert_2 \vert \vert \nabla b \vert \vert_2
\end{displaymath}
\end{theorem}
Here, $\vert \vert \cdot \vert \vert_{2,1}$ denotes the norm of the Lorentz space $L^{2,1}$. \\
A proof of this result can be found in \cite{CLMS} (see also \cite{W} and \cite{Tar}).\\

The next result - which can be found in \cite{Riv4} - gives an estimate of the Willmore energy of a compact surface with boundary. It is a consequence of Simon's monotonicity formula with boundary.

\begin{Lemma}\label{LemmaRiv1}
Let $\Sigma^2$ be a compact surface with boundary and let 
$\vec{\Phi}$ belong to $\mathcal{E}_{\Sigma^2}$. Then the
following inequality holds
\begin{equation}\label{LemmaTristan}
4 \pi \leq W(\vec{\Phi}) +2\frac{\mathcal{H}^1(\vec{\Phi}(\partial\Sigma^2))}{d(\vec{\Phi}(\partial \Sigma^2), \vec{\Phi}(\Sigma^2))}
\end{equation}
where $\mathcal{H}^1(\partial \vec{\Phi}(\Sigma^2))$ 
denotes the 1-dimensional Hausdorff measure of the boundary of the immersion $\vec{\Phi}(\partial\Sigma^2)$ and where $d( \cdot , \cdot )$ denotes the usual distance between two sets in $\mathbb{R}^3$.
\end{Lemma}
A useful lemma relating area, diameter and Willmore energy of a connected, compact surface without boundary is given by the following lemma due to Simon (see \cite{Si}). 

\begin{Lemma}\label{LemmaSimon}
Let $\Sigma^2\subset \R^3$  be an immersed connected, compact  surface without boundary.\\
Then it holds
\begin{displaymath}
\frac{Area(\Sigma^2)}{W(\Sigma^2)}\leq diam^2(\Sigma^2) \leq C Area(\Sigma^2)W(\Sigma^2)
\end{displaymath}
for some constant $C>0$.
\end{Lemma} 
Finally, let us recall the following lemma due to Simon (see \cite{Si}) giving an estimate for the Willmore energy from below. 

\begin{Lemma}\label{LemmaSimon2}
Assume that $\Sigma^2\subset \R^3$ is an immersed compact surface without boundary,  that $\partial B_{\rho}$ intersects $\Sigma^2$ transversely and that $\Sigma^2 \cap B_{\rho}$ contains disjoint subsets $\Sigma_1$, $\Sigma_2$ with $\Sigma_j \cap B_{\theta \rho} \neq \emptyset$, $\partial \Sigma_j \subset \partial B_{\rho}$, and $\vert \partial \Sigma_j \vert \leq \beta \rho$ for $j=1,2$, $\theta \in (0, \frac{1}{2})$ and $\beta > 0$. Then
\begin{displaymath}
W(\Sigma^2) \geq 8 \pi - C\beta \theta
\end{displaymath}
where $C$ does not depend on $\Sigma$, $\beta$ or $\theta$.
\end{Lemma}


\begin{thebibliography}{99}

\bibitem{A} David R. Adams, A note on Riesz potentials, Duke Math. J. 42 (1975), no. 4, 765-778 
\bibitem{AF} Robert A. Adams, John J. F. Fournier, Sobolev spaces, Elsevier, 2003
\bibitem{BK} Matthias Bauer, Ernst Kuwert, Existence of Minimizing Willmore Surfaces of Prescribed Genus, Int. Math. Res. Not. 10 (2003), 553-576  
\bibitem{BM} Giovanni Bellettini, Luca Mugnai, Approximation of Helfrich's functional via diffuse interfaces, SIAM J. Math. Anal. 42 (2010), no. 6, 2402-2433
\bibitem{BeRi} Yann Bernard, Tristan Rivi\`ere, Local Palais-Smale sequences for the Willmore functional, Comm. Anal. Geom. 19 (2011), no. 3, 563-599
\bibitem{BeRi2} Yann Bernard, Tristan Rivi\`ere, Energy Quantization for Willmore Surfaces and Applications, preprint (2011)
\bibitem{BeRi3} Yann Bernard, Tristan Rivi\`ere, Asymptotic Analysis of Branched Willmore Surfaces, to appear in Pacific J. Math. (2013)
\bibitem{Bla} Wilhelm Blaschke, ``Vorlesungen ¨\"uber Differentialgeometrie
und geometrische Grundlagen von Einsteins Relativit\"atstheorie'' III, Die Grundlehren der mathematischen
Wissenschaften in Einzeldarstellungen, Bd. XXIX Differentialgeometrie
der Kreise und Kugeln, bearbeitet von Gerhard
Thomsen, Springer, 1929
\bibitem{Bry} Robert Bryant, A duality theorem for Willmore surfaces, J. Differential Geom. 20 (1984), 23-53 
\bibitem{CL} Sagun Chanillo,  Yan Yan Li,  Continuity of solutions of uniformly elliptic equations in R2, Manuscripta Math. 77 (1992), no. 4, 415-433
\bibitem{C} Bang-Yen Chen, Some conformal Invariants of Submanifolds and Their Application, Boll. Unione Mat. Ital. 10 (1974), 380-385
\bibitem{C2} Bang-yen Chen,  On an inequality of T. J. Willmore,
Proc. Amer. Math. Soc. 26 1970 473-479.
\bibitem{CFS} Ralph Chill, Eva Fa\v sangov\'a, Reiner Sch\"atzle, Willmore blow-ups are never compact, Duke Math. J. 147 (2009), no. 2, 345-376
\bibitem{CV} Rustum Choksi, Marco Veroni, Global minimizers for the doubly-constrained Helfrich energy: the axisymmetric case, preprint (2012)
\bibitem{CLMS} R. Coifman, P. L. Lions, Y. Meyer, S. Semmes, Compensated compactness and Hardy spaces, J. Math. Pures Appl. 
72 (1993), 247-286
\bibitem{GT} David Gilbarg, Neil S. Trudinger, Elliptic Partial Differential Equations of Second Order, Springer, 1998
\bibitem{Hel} Fr\'ed\'eric H\'elein, Harmonic maps, conservation laws and moving frames, Cambridge Tracts in Mathematics, 2002
\bibitem{H} Wolfgang Helfrich, Elastic properties of lipid bilayers
- theory and possible experiments, Zeitschrift F\"ur Naturforschung
C - A Journal Of Biosciences. 28. (1973), 693-703
\bibitem{J} J\"urgen Jost, Compact Riemann Surfaces, Springer, 2006
\bibitem{K} Wilhelm P. A. Klingenberg, Riemannian Geometry, De Gruyter, 1995
\bibitem{KL} Ernst Kuwert, Yuxiang Li, $W^{2,2}$-conformal immersions of a closed Riemann Surface into $\R^n$, Comm. Anal. Geom. 20 (2012), no. 2, 313-340
\bibitem{KS1} Ernst Kuwert, Reiner Sch\"atzle, The Willmore flow with small initial energy, J. Differential Geom. 57 (2001), 409-441
\bibitem{KS2} Ernst Kuwert, Reiner Sch\"atzle, Gradient Flow for the Willmore Functional, Communications in Analysis and Geometry,
10 (2002), 307-339
\bibitem{KS3} Ernst Kuwert, Reiner Sch\"atzle, Removability of point singularities of Willmore surfaces, Ann. of Math. 160(2004), 315-357
\bibitem{KS4} Ernst Kuwert, Reiner Sch\"atzle, Branch points of Willmore surfaces, Duke Math. J. 138 (2007), no. 2, 179-201
\bibitem{KLS} Ernst Kuwert, Yuxiang Li, Reiner Sch\"atzle, The large genus limit of the infimum of the Willmore energy. Amer. J. Math. 132 (2010), no. 1, 37-51
\bibitem{KS5} Ernst Kuwert, Reiner Sch\"atzle, Closed surfaces with bound on their Willmore energy, to appear in Annali della Scuola Normale Superiore di Pisa, 2011, arXiv:math.DG/1009.5286.
\bibitem{KS6} Ernst Kuwert, Reiner Sch\"atzle, Minimizers of the Willmore functional under fixed conformal class, J. Differential Geom. 93 (2013), no. 3, 471-530
\bibitem{KS7} Ernst Kuwert, Reiner Sch\"atzle, The Willmore functional, Topics in modern regularity theory, 1--115, CRM Series, Scuola Normale Superiore Pisa, 2012
\bibitem{KMS} Ernst Kuwert, Andrea Mondino, Johannes Schygulla, Existence of immersed spheres minimizing curvature functionals in compact 3-manifolds, preprint, 2011
\bibitem{LMS} Tobias Lamm, Jan Metzger, Felix Schulze, Foliations of asymptotically flat manifolds by surfaces of Willmore type, Math. Ann. (2011), 350-378
\bibitem{LL} L. D. Landau, E. M. Lifshitz, Theory of elasticity
- Course of theoretical physics, volume 7 - third edition
revised and enlarged by E.M. Lifshitz, A.M. Kosevich and
L.P. Pitaevskii, Butterwoth-Heinemann 1986
\bibitem{LY} P. Li, S.-T.Yau, A New Conformal Invariant and its Applications
to the Willmore Conjecture and the First Eigenvalue on Compact Surfaces, Inventiones Math. 69 (1982), 269-291
\bibitem{LS} Reinhard Lipowsky, Erich Sackman, Structure and Dynamics of Membranes, Elsevier, 1995
\bibitem{MN} Fernando C. Marques, Andr\'e Neves, Min-Max theory and the Willmore conjecture, preprint (2012) 
\bibitem{Mil} Thomas Milcent, Shape derivative of the Willmore functional and
applications to equilibrium shapes of vesicles, Rapport de Recherche, 2011
\bibitem{MB1} X. Michalet, D. Bensimon, Fluctuating vesicles of nonspherical topology, Phys. Rev. Letter, 12(1), (1994), 168-171
\bibitem{MB2} X. Michalet, D. Bensimon, Observation of stable shapes and conformal diffusion in genus 2 vesicles, Science, 269,  (1995), 666-668
\bibitem{M} Andrea Mondino, Some results about the existence of critical points for the Willmore functional, Math. Z. 266 (2010), no. 3, 583-62
\bibitem{MR} Andrea Mondino, Tristan Rivi\`ere, Immersed Spheres of Finite Total Curvature into Manifolds, preprint (2012)
\bibitem{MR2} Andrea Mondino, Tristan Rivi\`ere, Willmore spheres in compact Riemannian manifolds, Advances Math. 232 (2013), no. 1, 608-676
\bibitem{MS} S. M\"uller, V. Sverak, On surfaces of finite total curvature, J. Differential Geom. 42 (1995), 229-258
\bibitem{MB} M. Mutz, D. Bensimon, Observation of toroidal vesicles, Phys. Rev. Letter A 43 (1991), 4525-4528 
\bibitem{PR} Mark A. Peletier, Matthias R\"oger, Partial Localization, Lipid Bilayers, and the Elastica Functional, Arch. Rational Mech. Anal. 193 (2009), 475-537
\bibitem{Mitoch} Christian Renken, Gino Siragusa, Guy Perkins, Lance Washington, Jim Nulton, Peter Salamon, Terrence G. Frey, A thermodynamic model describing the nature of the crista
junction: a structural motif in the mitochondrion, Journal of Structural Biology 138 (2002), 137-144\bibitem{Riv} Tristan Rivi\`ere, Conservation laws for conformal invariant variational problems, Invent. Math., 168 (2007), 1-22
\bibitem{Riv2} Tristan Rivi\`ere, Analysis aspects of Willmore surfaces, Inventiones Math., 174 (2008), no.1, 1-45 
\bibitem{Riv3}Tristan Rivi\`ere, Sequences of smooth global isothermic immersions, Comm. P. D. E. 38 (2013), no. 2, 276-303
\bibitem{Riv4} Tristan Rivi\`ere, Lipschitz conformal immersions from degenerating Riemann surfaces with $L^2$-bounded second fundamental forms, Adv. Calc. Var. 6 (2013), no. 1, 1-31
\bibitem{Riv5} Tristan Rivi\`ere, Variational Principles for immersed Surfaces with L2-bounded Second Fundamental Form, J. f\"ur die reine und angewandte Mathematik (Crelles journal) 2013
\bibitem{RivLect} Tristan Rivi\`ere, The role of conservation laws in the analysis of conformally invariant problems, in {\it Topics in Modern Regularity Theory}, G. Mingione Ed., CRM series 13, Edizioni Della Normale, 2012
\bibitem{RivNotes} Tristan Rivi\`ere, Conformally Invariant Variational Problems, in preparation 
\bibitem{RivOberw} Tristan Rivi\`ere, Error analysis for the Willmore-Helfrich Functional, Mini-workshop: Mathematics of Biological Membranes, Oberwolfach report 41/2008, 2305-2309
\bibitem{RivSt} Tristan Rivi\`ere, Michael Struwe, Partial regularity for harmonic maps, and relates problems, Comm. Pure Appl. Math.  61 (2008),  no. 4, 451-463
\bibitem{RS} Matthias R\"oger, Reiner Sch\"atzle, Control of the isoperimetric deficit by the Willmore deficit, Analysis (Munich) 32 (2012), no. 1, 1-7
\bibitem{RuSi} Thomas Runst, Winfried Sickel, Sobolev spaces of fractional order, Nemytskij operators, and nonlinear 
partial differential equations, De Gruyter, 1996
\bibitem{Sch1} Reiner Sch\"atzle, Lower semicontinuity of the Willmore functional for currents. J. Differential Geom. 81 (2009), no. 2, 437-456
\bibitem{Sch2} Reiner Sch\"atzle, The Willmore boundary problem. Calc. Var. Partial Differential Equations 37 (2010), no. 3-4, 275-302
\bibitem{Schy} Johannes Schygulla, Willmore minimizers with prescribed isoperimetric ratio, Arch. Ration. Mech. Anal. 203 (2012), no. 3, 901-941
\bibitem{Seifert} U. Seifert, R. Lipowsky, Morphology of vesicles, Handbook of Biological Physics, Elsevier, 1995
\bibitem{Rap} Yoko Shibata, Tom Shemesh, William A. Prinz, Alexander F. Palazzo, Michael M. Kozlov, Tom A. Rapoport, Mechanisms Determining the Morphology of the Peripheral ER, Cell 143 (2010), 775-788
\bibitem{Si} Leon Simon, Existence of surfaces minimizing the Willmore functional, Comm. Anal. Geom. 1 (1993), 281-326
\bibitem{S} Elias M. Stein, Singular Integrals and Differentiability Properties of Functions, Princeton University Press, 1970
\bibitem{S2} Elias M. Stein, Harmonic analysis, Princeton University Press, 1993
\bibitem{Tar} Luc Tartar, Remarks on oscillations and Stokes equation, in: "Macroscopic modeling of turbulent flows", Lect.
Notes in Physics 230, Springer 1985, 24-31
\bibitem{Ta} Michael E. Taylor, Partial Differential Equations I-III, Springer, 1996
\bibitem{T1} Peter Topping, Towards the Willmore conjecture, Calc. Var., 11 (2000), 361-393
\bibitem{T2} Peter Topping, An Approach to the Willmore Conjecture, in Global theory of minimal surfaces, ed. David Hoffman. Clay math. proc., vol. 2, AMS (2005), 769-77
\bibitem{U} Karen K. Uhlenbeck, Connections with $L^p$ bounds on curvature, Comm. Math. Phys. 83 (1982), no. 1, 31-42
\bibitem{Wei} Joel L. Weiner, On a problem of Chen, Willmore et al., Indiana University mathematics journal, 27 (1978), 19-35
\bibitem{W} Henry Wente, An existence theorem for surfaces of constant mean curvature, J. Math. Anal. Appl. 26
(1969), 318-344
\bibitem{Wil} Thomas J. Willmore, Riemannian Geometry, Oxford University Press, 1993
\bibitem{Wil2} Thomas K. Willmore, Total curvature in Riemannian geometry, 1982
\bibitem{ZH} Ou-Yang Zhong-can, Wolfgang Helfrich, Instability and Deformation of a Spherical Vesicle by Pressure, Phys. Rev. Lett, 59 no. 21, 1987
\end{thebibliography}
\end{document}